\documentclass[a4paper]{article}
\usepackage[english,russian]{babel}
\usepackage{latexsym,mathrsfs}
\usepackage{amssymb,amsmath,amsthm}
\usepackage{geometry,tikz}
\usepackage{stmaryrd}
\usepackage{bussproofs}
\usepackage{cancel}
\usepackage{color}
\usepackage{longtable}
\usepackage[urlbordercolor={.8 0 0}]{hyperref}
\usepackage{array}

\newcounter{rowcount}
\newcolumntype{T}{>{\small \ttfamily}l}
\newenvironment{hilbert}{ \setcounter{rowcount}{0}%
  \begin{center}
    \begin{tabular}{@{\stepcounter{rowcount} \makebox[3em][r]{\small \ttfamily \therowcount}\hspace*{\tabcolsep}}|l|T} 
 }
{\end{tabular} \end{center}}

\newtheorem{Prop}{Предложение}
\newtheorem{Lemm}[Prop]{Лемма}
\newtheorem{Th}[Prop]{Теорема}

\newtheorem{Cor}[Prop]{Следствие}

\theoremstyle{definition}

\newtheorem{Rem}[Prop]{Замечание}

\title{О кванторной версии модальной логики Белнапа--Данна и некоторых её расширениях}

\author{А.В.\ Грефенштейн}

\date{\normalsize Санкт-Петербургский государственный 
университет, Россия, 199034, Санкт-Петербург, Университетская наб. 7/9.\\ 
E-mail:\ aleksandrgrefenstejn@gmail.com;\ st082487@student.spbu.ru}

\begin{document}

\maketitle

\selectlanguage{english}
\begin{abstract}

We consider a quantified version of the (propositional) modal logic $\mathsf{BK}$, proposed earlier by S.\,P.\ Odintsov and H.\ Wansing; this version will be denoted by $\mathsf{QBK}$. Using the canonical model method, we prove the strong completeness of $\mathsf{QBK}$ with respect to a suitable possible world semantics with expanding domains. Similar results are obtained for some natural $\mathsf{QBK}$-extensions. In particular, it is proved that the extension of $\mathsf{QBK}$ with Barcan scheme is strongly complete with respect to a suitable possible world semantics with constant domains. Moreover, we define faithful embeddings (à la G\"{o}del--McKinsey--Tarski) of the quantified versions of Nelson's constructive logics into appropriate $\mathsf{QBK}$-extensions. 

\bigskip
{\footnotesize \emph{Keywords:} many-valued modal logic, strong negation, quatification, faithful embedding.
}
\end{abstract}
\selectlanguage{russian}

\tableofcontents

\def\fCenter{\rightarrow}

\vspace{4 mm}
\section{Введение}
В статье \cite{Odintsov-Wansing}  была рассмотрена пропозициональная логика $\mathsf{BK}$, представляющая собой белнапскую версию наименьшей нормальной модальной логики $\mathsf{K}$. Она представляет собой консервативное расширение $\mathsf{K}$ при помощи сильного, четырёхзначного отрицания $\sim$. С другой стороны, $\mathsf{BK}$ может быть получена путём добавления материальной импликации $(\rightarrow)$ и константы абсурд $(\bot)$ к наименьшей модальной логике $\mathsf{K}_{\mathsf{FDE}}$, обогащающей так называемое <<первоступенчатое следование>> ($\mathsf{FDE}$). С.\,П.\ Одинцовым и Х.\ Вансингом с помощью метода канонических моделей было доказано, что $\mathsf{BK}$ и некоторые её естественные расширения сильно полны относительно подходящей семантики Крипке. Установлено, что конструктивные логики Нельсона ${\mathsf{N4}}^{\bot}$ и $\mathsf{N3}$ точно вкладываются в подходящие расширения $\mathsf{BK}$ посредством аналога трансляции Маккинси--Гёделя--Тарского, погружающей пропозициональную интуиционистскую логику $\mathsf{Int}$ в модальную логику $\mathsf{S4}$.

Целью же данной работы является изучение кванторной версии $\mathsf{BK}$ и доказательство теоремы о сильной полноте для неё; обозначим эту логику  $\mathsf{QBK}$. При этом будет использоваться семантика типа Крипке с расширяющимися носителями. Сильная полнота также установлена для некоторых естественных расширений, полученных путём добавления новых схем аксиом или введения  ограничений на отношение достижимости в шкалах. В частности, разобран важный случай константных носителей. Как следствие теорем о полноте, предъявлено точное вложение кванторных логик Нельсона в подходящие расширения $\mathsf{QBK}$.
\section{Предварительные сведения}
Произвольное множество  не логических символов $\sigma$, каждому из которых поставлена в соответствие его арность, назовём сигнатурой. В рамках данной работы будем полагать, что в $\sigma$ нет функциональных символов. Обозначим за  $\mathrm{Pred}_{\sigma}$ и $\mathrm{Const}_{\sigma}$ соответственно множества предикатных и константных символов сигнатуры $\sigma$.
\medskip

Раз и навсегда  зафиксируем  счётное множество  \emph{переменных} $\mathrm{Var}$ и обозначим за  $\mathrm{Term}_{\sigma}$ множество $\sigma$-термов. Тогда мы получим, что $$\mathrm{Term}_{\sigma}=\mathrm{Var}\cup \mathrm{Const}_{\sigma}.$$. Будем использовать следующие логические символы:
\begin{itemize}

\item  \emph{символы пропозициональных связок} $\rightarrow$, $\wedge$, $\vee$, $\sim$;

\item  \emph{символы модальных операторов} $\square$ и $\lozenge$;

\item  \emph{символ логической константы} $\bot$;

\item  \emph{символы кванторов} $\forall$ и $\exists$.

\end{itemize}
Как всегда обозначим за $\mathrm{Form}_{\sigma}$  множество всех $\sigma$-формул  и  за $\mathrm{Sent}_{\sigma}$  множество всех $\sigma$-предложений, то есть формул, в которых отсутствуют свободные вхождения переменных. Произвольное множество $\sigma$-предложений $\Gamma$ будем называть \emph{$\sigma$-теорией}. Для удобства введём некоторые сокращения: $$\neg\Phi:=\Phi\rightarrow\bot,~~~~~~\Phi\leftrightarrow\Psi:=(\Phi\rightarrow\Psi)\wedge(\Psi\rightarrow\Phi),~~~~~~\Phi\Leftrightarrow\Psi:=(\Phi\leftrightarrow\Psi)\wedge(\sim\Phi\leftrightarrow{\sim\Psi}).$$

Произвольную функцию из множества $\mathrm{Var}$ в $\mathrm{Term}_\sigma$ будем называть $\sigma$-подстановкой. При этом $\sigma$-подстановку $\lambda$ называют \emph{основной}, если $\lambda \left( x \right) \in \mathrm{Const}_{\sigma}$
для любого $x \in \mathrm{Var}$. Пусть $\Phi$  такая $\sigma$-формула, что только переменные $x_1, \dots, x_n$ могут входить в неё свободно. В таком случае, если $\lambda$ --- это $\sigma$-подстановка, будем обозначать за
$\lambda \Phi$ результат одновременной подстановки в $\Phi$ $\sigma$-термов $\lambda \left( x_1 \right), \dots, \lambda \left( x_n \right)$ вместо всех свободных вхождений  $x_1, \ldots, x_n$ соответственно. В частности, если
\[
\lambda\ =\ {\left\{ \left( x, t \right) \right\} \cup \left\{ \left( y, y \right) \mid y \in Var ~ \text{и} ~ y \ne x \right\}}
\]
(где $x \in \mathrm{Var}$ и $t \in \mathrm{Term}_{\sigma}$), то пишем $\Phi \left( x/t \right)$ вместо $\lambda\Phi$.

Под $\sigma$-структурой будем подразумевать некоторое непустое множество с интерпретацией символов сигнатуры в нём. Пусть $\mathfrak{A}$ --- произвольная  $\sigma$-структура. Для любого $\varepsilon \in \sigma$ полагаем, что
\[
\varepsilon^{\mathfrak{A}}\ :=\ \text{интерпретация} ~ \varepsilon ~ \text{в} ~ \mathfrak{A} .
\]
 
В некоторых случаях будет удобно расширять $\sigma$ до
\[
\sigma_A\ :=\
{\sigma \cup \left\{ \underline{a} \mid a \in A \right\}},
\]
где $\underline{a}$ суть новые константы, и переходить от $\mathfrak{A}$ к её $\sigma_A$-расширению $\mathfrak{A}^{\ast}$, причём
\[
\underline{a}^{\mathfrak{A}^{\ast}}\ :=\ a \quad \text{для любых} \enskip {a \in A} .
\]
Будем называть $\sigma_A$-формулы просто \emph{$A$-формулами}. Далее, если $\Phi$ --- это $A$-предложение, то пишем
$\mathfrak{A} \models \Phi$ вместо $\mathfrak{A}^{\ast} \models \Phi$. Под \emph{$A$-подстановкой} подразумеваем $\sigma_A$-подстановку.

\section{Гильбертовское исчисление}

\medskip
Дедуктивная система нашей логики будет включать следующие аксиомы.

\begin{enumerate}
    \item  Аксиомы пропозициональной классической  логики  в языке $\{\vee,\wedge,\rightarrow,\bot\}$:

\begin{itemize}

\item[$\mathtt{I1}$.]
$\Phi \rightarrow \left( \Psi \rightarrow \Phi \right)$;

\item[$\mathtt{I2}$.]
$\left( \Phi \rightarrow \left( \Psi \rightarrow \Theta \right) \right) \rightarrow
\left( \left( \Phi \rightarrow \Psi \right) \rightarrow \left( \Phi \rightarrow \Theta \right) \right)$;

\item[$\mathtt{C1}$.]
$\Phi \wedge \Psi \rightarrow \Phi$;

\item[$\mathtt{C2}$.]
$\Phi \wedge \Psi \rightarrow \Psi$;

\item[$\mathtt{C3}$.]
$\Phi \rightarrow \left( \Psi \rightarrow {\Phi \wedge \Psi} \right)$;

\item[$\mathtt{D1}$.]
$\Phi \rightarrow \Phi \vee \Psi$;

\item[$\mathtt{D2}$.]
$\Psi \rightarrow \Phi \vee \Psi$;

\item[$\mathtt{D3}$.]
$\left( \Phi \rightarrow \Theta \right) \rightarrow
\left( \left( \Psi \rightarrow \Theta \right) \rightarrow \left( {\Phi \vee \Psi} \rightarrow \Theta \right) \right)$;

\item[$\mathtt{N1}$.]
$\Phi \vee \left(\Phi \rightarrow \bot\right)$;

\item[$\mathtt{N2}$.]
$\bot \rightarrow \Phi$.

\end{itemize}

\item Пропозициональные аксиомы сильного отрицания:

\begin{itemize}

\item[$\mathtt{SN1}$.]
$\sim\sim\Phi\leftrightarrow\Phi$;

\item[$\mathtt{SN2}$.] $\sim\left(\Phi\rightarrow\Psi\right)\leftrightarrow\left(\Phi\wedge\sim\Psi\right)$;

\item[$\mathtt{SN3}$.]$\sim\left(\Phi\vee\Psi\right)\leftrightarrow\left(\sim\Phi\wedge\sim\Psi\right)$;

\item[$\mathtt{SN4}$.]$\sim\left(\Phi\wedge\Psi\right)\leftrightarrow\left(\sim\Phi\vee\sim\Psi\right)$;

\item[$\mathtt{SN5}$.]$\sim\bot$.

\end{itemize}

\item Аксиомы $K$:

\begin{itemize}

\item[$\mathtt{K1}$.]$\left(\square\Phi\wedge\square\Psi\right)\rightarrow\square\left(\Phi\wedge\Psi\right)$; 

\item[$\mathtt{K2}$.]$\square\left(\Phi\rightarrow\Phi\right)$.
    
\end{itemize}

\item Аксиомы связи модальностей:

\begin{itemize}

\item[$\mathtt{M1}$.]$\neg\square\Phi\leftrightarrow\lozenge\neg\Phi$;

\item[$\mathtt{M2}$.]$\neg\lozenge\Phi\leftrightarrow\square\neg\Phi$;

\item[$\mathtt{M3}$.]${\square\Phi}\Leftrightarrow{\sim\lozenge\sim\Phi}$;

\item[$\mathtt{M4}$.]$\lozenge\Phi\Leftrightarrow{\sim\square\sim\Phi}$.
    
\end{itemize}

\item Кванторные аксиомы:

\begin{itemize}

\item[$\mathtt{Q1}$.]$\forall x\,\Phi\rightarrow\Phi(x/t)$, \  где $t$  свободен для $x$ в $\Phi$;

\item[$\mathtt{Q2}$.]$\Phi(x/t)\rightarrow\exists x\,\Phi$, \  где $t$  свободен для $x$ в $\Phi$;

\item[$\mathtt{Q3}$.] $\sim\forall x\, \Phi\leftrightarrow\exists x\sim\Phi$;

\item[$\mathtt{Q4}$.] $\sim\exists x\, \Phi\leftrightarrow\forall x\sim\Phi$.

\end{itemize}
\end{enumerate}

Также будем использовать следующие правила вывода:
\begin{itemize}

\item[$\mathtt{MP}$.] \emph{modus ponens}\
\begin{prooftree}
\AxiomC{$\Phi$}
\AxiomC{$\Phi \rightarrow \Psi$}
\RightLabel{;}
\BinaryInfC{$\Psi$}
\end{prooftree}

\item[$\mathtt{MB}$.]  \emph{правило монотонности для $\square$  }
\begin{prooftree}
\Axiom$\Phi \fCenter \Psi$
\RightLabel{;}
\UnaryInf$\square\Phi \fCenter {{\square} \Psi}$
\end{prooftree}

\item[$\mathtt{MD}$.]  \emph{правило монотонности для $\lozenge$  }
\begin{prooftree}
\Axiom$\Phi \fCenter \Psi$
\RightLabel{;}
\UnaryInf$\lozenge\Phi \fCenter {{\lozenge} \Psi}$
\end{prooftree}

\item[$\mathtt{BR1}$.]  \emph{правило Бернайса для  $\forall$}
\begin{prooftree}
\Axiom$\Phi \fCenter \Psi$
\RightLabel{,\quad где переменная $x$ не свободна в $\Phi$;}
\UnaryInf$\Phi \fCenter {{\forall x}\, \Psi}$
\end{prooftree}

\item[$\mathtt{BR2}$.]  \emph{правило Бернайса для  $\exists$}
\begin{prooftree}
\Axiom$\Phi \fCenter \Psi$
\RightLabel{,\quad где переменная $x$ не свободна в $\Psi$.}
\UnaryInf${{\exists x}\, \Phi} \fCenter \Psi$
\end{prooftree}

\end{itemize}

Тогда $\mathsf{QBK}_\sigma$ --- наименьшее множество $\sigma$-формул, содержащее перечисленные выше схемы аксиом и замкнутое относительно правил вывода $\mathtt{MP}$, $\mathtt{MB}$, $\mathtt{MD}$, $\mathtt{BR1}$ и $\mathtt{BR2}$\footnote{Всякий раз, когда будет понятно, о какой сигнатуре идёт речь, или когда это просто будет удобно, будем опускать нижний индекс $\sigma$.}.

\medskip
Для любых $\Gamma \subseteq \mathrm{Form}_{\sigma}$ определим
\[
{\mathrm{Disj} \left( \Gamma \right)}\ :=\
{\left\{
\Phi_1 \vee \dots \vee \Phi_n \mid \left\{ \Phi_1, \dots, \Phi_n \right\} \subseteq \Gamma,~n\in\mathbb{N}
\right\}} .
\]
При этом пустую дизъюнкцию отождествим с $\bot$.

\medskip
Пусть $\Gamma \subseteq \mathrm{Sent}_{\sigma}$ и $\Delta \subseteq \mathrm{Form}_{\sigma}$, будем писать $\Gamma \vdash \Delta$ если и только если некоторый элемент
$\mathrm{Disj} \left( \Delta \right)$ может быть получен из элементов  $\Gamma \cup \mathsf{QBK}_{\sigma}$ с помощью  $\mathtt{MP}$, $\mathtt{BR1}$
и $\mathtt{BR2}$.

\vspace{6 mm}

\begin{Th}[о дедукции] \label{Deduction_Theorem}
Для любых $\Gamma \cup \left\{ \Phi \right\} \subseteq \mathrm{Sent}_{\sigma}$ и $\Psi
\in \mathrm{Form}_{\sigma}$,
\[
{\Gamma \cup \left\{ \Phi \right\} \vdash \Psi}
\quad \Longleftrightarrow \quad
{\Gamma \vdash \Phi \rightarrow \Psi} .
\]
\end{Th}

\begin{proof}
Доказательство аналогично случаю логики предикатов.
\end{proof}

\vspace{4 mm}

Отметим также один достаточно очевидный факт:




\begin{Prop}[нормализация] Правило

\begin{prooftree}
\AxiomC{$\Phi$}
\RightLabel{$(\mathsf{RN})$}
\UnaryInfC{$\square\Phi$}
\end{prooftree}
является допустимым в $\mathsf{QBK}_\sigma$.
\begin{proof}

Пусть $\Phi\in \mathsf{QBK}_{\sigma}$. Тогда по аксиоме $\mathtt{I1}$ и правилу $\mathtt{MP}$ получаем, что формула  $(\Phi\rightarrow\Phi)\rightarrow\Phi$ лежит в $\mathsf{QBK}_{\sigma}$. По аксиоме  $\mathtt{K2}$ получаем также, что $\square(\Phi\rightarrow\Phi)$ лежит в $\mathsf{QBK}_{\sigma}$. По правилу монотонности для $\square$ и $\mathtt{MP}$ имеем: $\square\Phi\in \mathsf{QBK}_{\sigma}.$

\end{proof}

\end{Prop}

\vspace{4 mm}

Пусть $\Phi,\Psi,\Theta\in \mathrm{Form}_\sigma$. Тогда результат одновременной  замены всех вхождений формулы $\Phi$ в $\Theta$ на $\Psi$ обозначим за $\Theta(\Phi/\Psi)$.

\begin{Prop}[позитивное правило замены]
Правило 

\begin{prooftree}
\AxiomC{$\Phi\ \leftrightarrow\ \Psi$}
\RightLabel{$(\mathsf{PR})$,}
\UnaryInfC{$ \Theta\ \leftrightarrow \ \Theta(\Phi/\Psi)$}
\end{prooftree}
где $\Theta$  не содержит $\sim$, является допустимым в $\mathsf{QBK}_\sigma$.

\end{Prop}

\begin{proof}
Индукция по сложности формулы $\Theta$:

\medskip
1) Для атомарных формул результат очевиден ввиду единственной подформулы.

\medskip
2) Так как $\mathsf{QBK}_\sigma$ содержит аксиомы классической пропозициональной логики, то формулы $$(\Theta_1 \leftrightarrow\ \Theta_1(\Phi/\Psi))\rightarrow((\Theta_2 \leftrightarrow\ \Theta_2(\Phi/\Psi))\rightarrow((\Theta_1 * \Theta_2)\leftrightarrow(\Theta_1(\Phi/\Psi) * \Theta_2(\Phi/\Psi)))),   $$
где $*\in\{\vee,\wedge,\rightarrow\}$ лежат в $\mathsf{QBK}_\sigma$. Поэтому случай $\Theta=\Theta_1*\Theta_2$ получается тривиально.

\medskip
3) Случаи $\Theta=*\Theta_1$, где $*\in\{\square,\lozenge\}$, получаются применением соответствующего правила монотонности.

4) Пусть $\Theta=\exists x \Theta_1$. Тогда, за исключением тривиальной замены,  $\Theta(\Phi/\Psi)=\exists x \Theta_1(\Phi/\Psi)$. По предположению индукции $\Theta_1\leftrightarrow\Theta_1(\Phi/\Psi)$. Покажем, что $\exists x\Theta_1\rightarrow\exists x\Theta_1(\Phi/\Psi)$ выводима

\begin{hilbert}

$\Theta_1\rightarrow\Theta_1(\Phi/\Psi)$& индукционное предположение \\
$\Theta_1(\Phi/\Psi)\rightarrow\exists x\Theta_1(\Phi/\Psi)$& $\mathtt{Q2}$\\
$\Theta_1\rightarrow\exists x\Theta_1(\Phi/\Psi)$& транзитивность + 1 + 2\\
$\exists x\Theta_1\rightarrow\exists x\Theta_1(\Phi/\Psi)$& $\mathtt{BR2}$

\end{hilbert}
В обратную сторону аналогично.

5) Пусть $\Theta=\forall x \Theta_1$. Тогда, за исключением тривиальной замены,   $\Theta(\Phi/\Psi)=\forall x \Theta_1(\Phi/\Psi)$. Покажем, что $\forall x\Theta_1\rightarrow\forall x\Theta_1(\Phi/\Psi)$ выводима

\begin{hilbert}
$\forall x \Theta_1\rightarrow\Theta_1(\Phi/\Psi)$& $\mathtt{Q1}$\\
$\Theta_1\rightarrow\Theta_1(\Phi/\Psi) $& индукционное предположение\\
$\forall x \Theta_1\rightarrow\Theta_1(\Phi/\Psi)$& транзитивность + 1 + 2\\
$\forall x \Theta_1\rightarrow\forall x\Theta_1(\Phi/\Psi)$& $\mathtt{BR1}$
\end{hilbert}
В обратную сторону аналогично.

\end{proof}
\begin{Prop}[слабое правило замены]
Правило

\begin{prooftree}
\AxiomC{$\Phi \ \Leftrightarrow \ \Psi $}
\RightLabel{$(\mathsf{WR})$}
\UnaryInfC{$\Theta\ \Leftrightarrow \  \Theta(\Phi/\Psi)$}
\end{prooftree}
является допустимым в $\mathsf{QBK}_\sigma$.

\end{Prop}

\begin{proof} Достаточно доказать ослабленный вариант правила, а именно: \begin{prooftree}
\AxiomC{$\Phi \ \Leftrightarrow \ \Psi $}
\RightLabel{$(\widetilde{\mathsf{WR}}).$}
\UnaryInfC{$\Theta\ \leftrightarrow \  \Theta(\Phi/\Psi)$}
\end{prooftree} Действительно, чтобы получить $\sim\Theta\ \leftrightarrow \  \sim\Theta(\Phi/\Psi)$, нужно применить $\widetilde{\mathsf{WR}}$ к формуле $\sim\Theta$. \vspace{2 mm} \\ Сразу отметим, что если заменяется сама $\Theta$, то результат следует из $\Phi\Leftrightarrow\Psi$, поэтому далее такие случаи рассматривать не будем.

\medskip
Индукция по сложности формулы $\Theta$. \\
База очевидна. Ясно, что если $\Theta$ не начинается с $\sim$, то шаг индукции проводится аналогично индукционному переходу в доказательстве позитивного правила замены. Поэтому далее рассматриваем только $\Theta$, начинающиеся с $\sim$.

\medskip
1) Если $\Theta={\sim P(t_1,\ldots,t_n)}={\sim\Phi}$, то $\Theta(\Phi/\Psi)={\sim\Psi}$. И так как $\sim\Phi\leftrightarrow{\sim\Psi}$,  имеем $\Theta\leftrightarrow\Theta(\Phi/\Psi)$. 
\medskip

2) Если $\Theta={\sim(\Theta_1*\Theta_2)}$, где $*\in\{\wedge,\vee,\rightarrow\}$, то для каждой связки результат получается почти идентичным образом. Для примера разберём случай $\Theta={\sim(\Theta_1\wedge\Theta_2)}$: \vspace{2 mm} \\
Если заменяется $\Theta_1\wedge\Theta_2$, то все получается идентично случаю 1). В противном случае имеем $\Theta(\Phi/\Psi)=\sim(\Theta_1(\Phi/\Psi)\wedge\Theta_2(\Phi/\Psi))$. По аксиоме $\mathtt{SN4}$ $$\sim(\Theta_1(\Phi/\Psi)\wedge\Theta_2(\Phi/\Psi))\leftrightarrow(\sim\Theta_1(\Phi/\Psi)\vee\sim\Theta_2(\Phi/\Psi)),~~~~~\sim(\Theta_1\wedge\Theta_2)\leftrightarrow(\sim\Theta_1\vee\sim\Theta_2).$$
И по предположению индукции $$\sim\Theta_1\leftrightarrow\sim\Theta_1(\Phi/\Psi),~~~~~\sim\Theta_2\leftrightarrow\sim\Theta_2(\Phi/\Psi).$$
Значит, используя те же тавтологии классической логики как в  доказательстве $(\mathsf{PR})$, легко получаем, что $$(\sim\Theta_1\vee\sim\Theta_2)\leftrightarrow(\sim\Theta_1(\Phi/\Psi)\vee\sim\Theta_2(\Phi/\Psi)).$$
А значит и $$\Theta={\sim{(\Theta_1\wedge\Theta_2)}}\leftrightarrow{\sim(\Theta_1(\Phi/\Psi)\wedge\Theta_2(\Phi/\Psi))}=\Theta(\Phi/\Psi).$$
Для других случаев проводим аналогичное рассуждение, используя аксиомы $\mathtt{SN2}$ и $\mathtt{SN3}$.

\medskip
3) Пусть $\Theta=\sim*\Theta_1$, где $*\in\{\square,\lozenge\}$.  Для примера рассмотрим случай с $\lozenge$: \vspace{2 mm} \\
Если заменяется $\lozenge\Theta_1$, то тривиально.
Иначе $\Theta(\Phi/\Psi)=\sim\lozenge\Theta_1(\Phi/\Psi)$. По аксиоме $\mathtt{M4}$ имеем $$\sim\lozenge\Theta_1(\Phi/\Psi)\leftrightarrow\sim\sim\square\sim\Theta_1(\Phi/\Psi),~~~~~\sim\lozenge\Theta_1\leftrightarrow\sim\sim\square\sim\Theta_1.$$
Далее, по аксиоме $\mathtt{SN1}$ $$\sim\sim\square\sim\Theta_1(\Phi/\Psi)\leftrightarrow\square\sim\Theta_1(\Phi/\Psi),~~~~~\sim\sim\square\sim\Theta_1\leftrightarrow\square\sim\Theta_1.$$
И по предположению индукции $$\sim\Theta_1(\Phi/\Psi)\leftrightarrow\sim\Theta_1.$$
Тогда всё получается по правилу монотонности для $\square$. В случае с $\square$ действуем аналогично, используя аксиомы $\mathtt{M3}$, $\mathtt{SN1}$ и правило монотонности для $\lozenge$.

4) Пусть $\Theta={\sim*x\Theta_1}$, где $*\in\{\exists,\forall\}$.  Для примера рассмотрим случай с $\exists$: \vspace{2 mm} \\
Если заменяется $\forall x\,\Theta_1$, то снова тривиально. Иначе  $\Theta(\Phi/\Psi)={\sim\exists x\,\Theta_1(\Phi/\Psi)}$. По аксиоме $\mathtt{Q4}$ получаем $$\sim\exists x\,\Theta_1(\Phi/\Psi)\leftrightarrow\forall x\sim\Theta_1(\Phi/\Psi),~~~~~\sim\exists x\,\Theta_1\leftrightarrow\forall x\sim\Theta_1.$$
И по предположению индукции $$\sim\Theta_1(\Phi/\Psi)\leftrightarrow{\sim\Theta_1}.$$ 
Теперь, повторив рассуждения в пункте 5) из доказательства  допустимости позитивного правила замены, имеем $$\forall x\sim\Theta_1(\Phi/\Psi)\leftrightarrow\forall x\sim\Theta_1.$$ 
Таким образом, все случаи разобраны\footnote{При использовании правил $(\mathsf{PR})$ и $(\mathsf{WR})$, вместо замены всех вхождений сразу можно заменять какие-то конкретные. Далее будем пользоваться этой возможностью. }.
\end{proof}

\begin{Prop}
\ $\neg\neg\Phi\Leftrightarrow\sim\neg\Phi$

\end{Prop}

\begin{proof}

Сначала докажем, что $\vdash\neg\neg\Phi\rightarrow\sim\neg\Phi$

\medskip
\begin{hilbert}

$\neg\neg\Phi\rightarrow\Phi$& классическая тавтология \\
$\Phi\rightarrow(\sim\bot\rightarrow(\Phi\wedge\sim\bot))$& $\mathtt{C3}$\\
$\neg\neg\Phi\rightarrow(\sim\bot\rightarrow(\Phi\wedge\sim\bot))$& транзитивность  1 и 2\\
$\sim\bot\rightarrow(\neg\neg\Phi\rightarrow(\Phi\wedge\sim\bot))$& перестановка посылок в 3\\
$\sim\bot$& $\mathtt{SN5}$\\
$\neg\neg\Phi\rightarrow(\Phi\wedge\sim\bot)$& $\mathtt{MP}$ 4 и 5\\
$(\Phi\wedge\sim\bot)\rightarrow\sim(\Phi\rightarrow\bot)$& $\mathtt{SN2}$\\
$\neg\neg\Phi\rightarrow\sim(\Phi\rightarrow\bot)$& транзитивность  6 и 7

\end{hilbert}
И так как $$\neg\neg\Phi\rightarrow\sim(\Phi\rightarrow\bot)=\neg\neg\Phi\rightarrow\sim\neg\Phi,$$
мы получили вывод в одну сторону . \vspace{2 mm} \\ 
Теперь докажем, что $\vdash{\sim\neg\Phi}\rightarrow\neg\neg\Phi$

\begin{hilbert}

$\sim(\Phi\rightarrow\bot)\rightarrow(\Phi\wedge\sim\bot)$& $\mathtt{SN2}$ \\
$(\Phi\wedge\sim\bot)\rightarrow\Phi$& $\mathtt{C1}$\\
$\sim(\Phi\rightarrow\bot)\rightarrow\Phi$& транзитивность 1 и 2\\
$\Phi\rightarrow\neg\neg\Phi$& классическая тавтология\\
$\sim(\Phi\rightarrow\bot)\rightarrow\neg\neg\Phi$& транзитивность 3 и 4

\end{hilbert}
И так как $$\sim(\Phi\rightarrow\bot)\rightarrow\neg\neg\Phi=\sim\neg\Phi\rightarrow\neg\neg\Phi,$$
мы получили вывод в другую сторону.\vspace{2 mm} \\ 
Теперь заметим, что $$\sim\sim\neg\Phi\leftrightarrow\neg\Phi\leftrightarrow\neg\neg\neg\Phi\leftrightarrow\sim\neg\neg\Phi.$$
Таким образом, $$\neg\neg\Phi\Leftrightarrow\sim\neg\Phi.$$

\end{proof}

\begin{Cor}[сильная эквивалентность некоторых аксиом] 

Формулы, стоящие в левой и правой частях аксиом $\mathtt{SN1}$,  $\mathtt{SN3}$,  $\mathtt{SN4}$,  $\mathtt{Q3}$, $\mathtt{Q4}$, $\mathtt{M1}$ и $\mathtt{M2}$, сильно эквивалентны.

\end{Cor}

\begin{proof} Случай $\mathtt{SN1}$ получается задаром так как если хотим получить $$\sim\sim\sim\Phi\leftrightarrow\sim\Phi,$$
то достаточно подставить в $\mathtt{SN1}$ $\Psi=\sim\Phi$.

\medskip
$\mathtt{SN3}$: Нужно доказать $\sim\sim(\Phi\vee\Psi)\leftrightarrow\sim(\sim\Phi\wedge\sim\Psi)$. Имеем $$(\Phi\vee\Psi)\leftrightarrow\sim\sim(\Phi\vee\Psi),~~~~~\sim(\sim\Phi\wedge\sim\Psi)\leftrightarrow(\sim\sim\Phi\vee\sim\sim\Psi)\leftrightarrow(\Phi\vee\Psi).$$

\medskip
$\mathtt{SN4}$: Нужно доказать  $\sim\sim(\Phi\wedge\Psi)\leftrightarrow\sim(\sim\Phi\vee\sim\Psi)$. Имеем $$(\Phi\wedge\Psi)\leftrightarrow\sim\sim(\Phi\wedge\Psi),~~~~~\sim(\sim\Phi\vee\sim\Psi)\leftrightarrow(\sim\sim\Phi\wedge\sim\sim\Psi)\leftrightarrow(\Phi\wedge\Psi).$$

\medskip
$\mathtt{Q3}$: Нужно доказать $\sim\sim\forall x \Phi\leftrightarrow\sim\exists x\sim\Phi$. Имеем $$\forall x\Phi\leftrightarrow\sim\sim\forall x \Phi,~~~~~\sim\exists x\sim\Phi\leftrightarrow\forall x \sim\sim\Phi\leftrightarrow\forall x \Phi .$$

\medskip
$\mathtt{Q4}$: Нужно доказать $\sim\sim\exists x \Phi\leftrightarrow\sim\forall x\sim\Phi$. Имеем $$\exists x\Phi\leftrightarrow\sim\sim\exists x \Phi,~~~~~\sim\forall x\sim\Phi\leftrightarrow\exists x \sim\sim\Phi\leftrightarrow\exists x \Phi .$$

\medskip
$\mathtt{M1}$: Нужно доказать $\sim\neg\square\Phi\leftrightarrow\sim\lozenge\neg\Phi$. Имеем $$\neg\neg\square\Phi\leftrightarrow\square\Phi.$$
С помощью закона контрапозиции, снятия двойного отрицания и аксиомы $\mathtt{M2}$ получаем, что $$\square\Phi\leftrightarrow\neg\lozenge\neg\Phi\leftrightarrow\square\neg\neg\Phi.$$
Значит, $$\neg\neg\square\Phi\leftrightarrow\square\neg\neg\Phi.$$
Так как  $\neg\neg\Phi\Leftrightarrow\sim\neg\Phi$, то, применяя слабое правило замены, имеем $$\sim\neg\square\Phi\leftrightarrow\square\sim\neg\Phi.$$
По аксиоме $\mathtt{M4}$ и $\mathtt{SN1}$ получаем, что $$\sim\neg\square\Phi\leftrightarrow\sim\lozenge\neg\Phi.$$
Случай аксиомы $\mathtt{M2}$ получается аналогично.

\end{proof}

\begin{Prop} $(\Phi\rightarrow\Psi)\Leftrightarrow(\neg\Phi\vee\Psi)$

\end{Prop}

\begin{proof}
$(\Phi\rightarrow\Psi)\leftrightarrow(\neg\Phi\vee\Psi)$ есть как классическая тавтология. Разберёмся с сильным отрицанием. По аксиоме $\mathtt{SN2}$\, $\sim(\Phi\rightarrow\Psi)\leftrightarrow(\Phi\wedge\sim\Psi)$. По слабому правилу замены и предложению 6 получаем, что $(\neg\neg\Phi\wedge\sim\Psi)\leftrightarrow(\sim\neg\Phi\wedge\sim\Psi).$ И по аксиоме $\mathtt{SN3}$\, $(\sim\neg\Phi\wedge\sim\Psi)\leftrightarrow{\sim(\neg\Phi\vee\Psi)}.$  Тогда достаточно доказать, что $(\Phi\wedge\sim\Psi)\leftrightarrow(\neg\neg\Phi\wedge\sim\Psi)$. Но это частный случай классической тавтологии.
\end{proof}

Говорим, что формула $\Phi$ находится в негативной нормальной форме, если символ $\sim$ стоит только перед атомарными подформулами или перед константой $\bot$.

\begin{Th}[о негативной нормальной форме] Для любой формулы $\Phi\in \mathrm{Form}_{\sigma}$ существует формула $\overline{\Phi}\in \mathrm{Form}_{\sigma}$ в негативной нормальной форме такая, что $$\Phi\Leftrightarrow\overline{\Phi}\in \mathsf{QBK}_{\sigma}.$$

\end{Th}

\begin{proof}
Проводим индукцию по построению формулы $\Phi$, используя установленные сильные эквивалентности аксиом, предложения 5 и 7 и слабое правило замены.

\end{proof}
\section{Семантическое следование}
 Под  \emph{шкалой} будем подразумевать пару $\mathcal{W}=\langle W,R\rangle$, где $W$ --- непустое множество <<миров>> и $R\subseteq W\times W$. \vspace{2 mm} \\
Каждому миру произвольной  шкалы $\mathcal{W}$ поставим в соответствие две $\sigma$-структуры $\mathfrak{A}_w^+$ и $\mathfrak{A}_w^-$ такие, что $A_w^+=A_w^-=A_w\ $ и $~~c^{\mathfrak{A}_w^+}=c^{\mathfrak{A}_w^-}$ для любой константы $c\in \mathrm{Const}_{\sigma}$. Обозначим также $\mathscr{A}^+ = \langle \mathfrak{A}_w^+ \mid w \in W 
\rangle$ и $\mathscr{A}^- = \langle \mathfrak{A}_w^- \mid w \in W 
\rangle$. $\mathscr{A}^+$ и $\mathscr{A}^-$ можно рассматривать как некоторые функции из множества миров в класс всех $\sigma$-структур.\vspace{2 mm} \\
Упорядоченная тройка $\mathcal{M}=\langle\mathcal{W},\mathscr{A}^+,\mathscr{A}^-\rangle$ называется \emph{$\mathsf{QBK}_{\sigma}$-моделью}, если для любых $ u,v\in W$

\begin{itemize}
    \item $uRv\ \implies \  A_u\subseteq A_v$; 
    \item $uRv \ \implies \ \text{для любой}~c\in \mathrm{Const}_{\sigma}~~c^{\mathfrak{A}_u^+}=c^{\mathfrak{A}_v^+}$.
\end{itemize}

\medskip

Определим отношения $\Vdash^+$ и $\Vdash^-$ между парой (модель,мир) и $A_w$-предложением $\Theta$ \\ ($\mathcal{M},w\Vdash^+\Theta$, $\mathcal{M},w\Vdash^-\Theta$)  индукцией по построению  $\Theta$:

\medskip
\begin{center}
Если $\Theta$ атомарное, то $\mathcal{M},w\Vdash^+\Theta \iff \mathfrak{A}^+_{w}\Vdash\Theta$~ и ~$\mathcal{M},w\Vdash^-\Theta \iff \mathfrak{A}^-_{w}\Vdash\Theta$.
\end{center}
\begin{align*}
&\mathcal{M},w\Vdash^+\Phi\wedge\Psi\iff \mathcal{M},w\Vdash^+\Phi ~\text{и}~ \mathcal{M},w\Vdash^+\Psi; \\
&\mathcal{M},w\Vdash^-\Phi\wedge\Psi\iff \mathcal{M},w\Vdash^-\Phi ~\text{или}~ \mathcal{M},w\Vdash^-\Psi; \\
\\
&\mathcal{M},w\Vdash^+\Phi\vee\Psi\iff\ \mathcal{M},w\Vdash^+\Phi ~\text{или}~ \mathcal{M},w\Vdash^+\Psi; \\
&\mathcal{M},w\Vdash^-\Phi\vee\Psi\iff\ \mathcal{M},w\Vdash^-\Phi~ \text{и}~ \mathcal{M},w\Vdash^-\Psi; \\
\\
&\mathcal{M},w\Vdash^+\Phi\rightarrow\Psi\iff\ \mathcal{M},w\nVdash^+\Phi~\text{или}~\mathcal{M},w\Vdash^+\Psi; \\
&\mathcal{M},w\Vdash^-\Phi\rightarrow\Psi\iff\ \mathcal{M},w\Vdash^+\Phi ~\text{и}~ \mathcal{M},w\Vdash^-\Psi; \\
\\
&\mathcal{M},w\Vdash^+\square\Phi\iff\forall u\in W~(wRu\implies\mathcal{M},w\Vdash^+\Phi); \\
&\mathcal{M},w\Vdash^-\square\Phi\iff\exists  u\in W~(wRu ~\text{и}~\mathcal{M},w\Vdash^-\Phi); \\
\\
&\mathcal{M},w\Vdash^+\lozenge\Phi\iff\exists u\in W~(wRu ~\text{и}~\mathcal{M},w\Vdash^+\Phi); \\
&\mathcal{M},w\Vdash^-\lozenge\Phi\iff\forall u\in W~(wRu\implies\mathcal{M},w\Vdash^-\Phi); \\
\\
&\mathcal{M},w\Vdash^+\sim\Phi\iff\mathcal{M},w\Vdash^-\Phi; \\
&\mathcal{M},w\Vdash^-\sim\Phi\iff\mathcal{M},w\Vdash^+\Phi; \\
\\
&\mathcal{M},w\Vdash^+\forall x\Phi\iff\forall a\in A_w~ \mathcal{M},w\Vdash^+\Phi(x/\underline{a}); \\
&\mathcal{M},w\Vdash^-\forall x\Phi\iff\exists a\in A_w~ \mathcal{M},w\Vdash^-\Phi(x/\underline{a}); \\
\\
&\mathcal{M},w\Vdash^+\exists x\Phi\iff\exists a\in A_w~ \mathcal{M},w\Vdash^+\Phi(x/\underline{a}); \\
&\mathcal{M},w\Vdash^-\exists x\Phi\iff\forall a\in A_w~ \mathcal{M},w\Vdash^-\Phi(x/\underline{a}). \\
\end{align*}

\begin{center}
Также для $\bot$ полагаем:  $\mathcal{M},w\Vdash^-\bot$; \ \ \ $\mathcal{M},w\nVdash^+\bot$.
\end{center}
Будем называть отношения $\Vdash^+$ и\, $\Vdash^-$ отношением \emph{верифицируемости} и \emph{фальсифицируемости} соответственно. Тогда запись $\mathcal{M},w\Vdash^+\Theta$ ($\mathcal{M},w\Vdash^-\Theta$) будет читаться так: <<предложение $\Theta$ верифицируемо (фальсифицируемо) в мире $w$ модели $\mathcal{M}$>>.   
\medskip

Пусть $\Gamma\subseteq \mathrm{Sent}_{\sigma}$ и $\Delta\subseteq \mathrm{Form}_{\sigma}$. Говорят, что $\Delta$ \emph{семантически следует} из $\Gamma$ ($\Gamma\vDash\Delta$),  если и только если для любой $\mathsf{QBK}_{\sigma}$-модели $\mathcal{M}=\langle\mathcal{W},\mathscr{A}^+,\mathscr{A}^-\rangle$, мира $w\in W$  и основной  $A_w$-подстановки $\lambda$
\[
{\mathcal{M}, w \Vdash^+ \Phi} \quad \text{для любых} \enskip {\Phi \in \Gamma}
\quad \Longrightarrow \quad
{\mathcal{M}, w \Vdash^+ \lambda \Psi} \quad \text{для некоторой} \enskip {\Psi \in \Delta} .
\]

\begin{Th}[О семантической дедукции]
Для любых $\Gamma \cup \left\{ \Phi \right\} \subseteq \mathrm{Sent}_{\sigma}$ и $\Psi
\in \mathrm{Form}_{\sigma}$,
\[
{\Gamma \cup \left\{ \Phi \right\} \vDash \Psi}
\quad \Longleftrightarrow \quad
{\Gamma \vDash \Phi \rightarrow \Psi} .
\]
\end{Th}

\begin{proof}
Пусть $\mathcal{M}$ --- произвольная $\mathsf{QBK}_{\sigma}$-модель, $w\in W$  и $\lambda$ --- основная $A_w$-подстановка. Тогда, расписывая по определению правую часть утверждения, имеем $$\mathcal{M},w\Vdash^+\Theta ~~\text{для любых}~~\Theta\in\Gamma\implies\mathcal{M},w\Vdash^+\lambda(\Phi\rightarrow\Psi)=\Phi\rightarrow\lambda(\Psi).$$ Это эквивалентно тому, что$$\mathcal{M},w\Vdash^+\Theta ~~\text{для любых}~~\Theta\in\Gamma\implies\left(\mathcal{M},w\Vdash^+\Phi\implies\mathcal{M},w\Vdash^+\lambda(\Psi)\right).$$ А это эквивалентно тому, что $$\mathcal{M},w\Vdash^+\Theta ~~\text{для любых}~~\Theta\in\Gamma\cup\{\Phi\}\implies\mathcal{M},w\Vdash^+\lambda(\Psi).$$
Поскольку модель, мир и подстановка были произвольны, получаем, что \[
{\Gamma \cup \left\{ \Phi \right\} \vDash \Psi}
\quad \Longleftrightarrow \quad
{\Gamma \vDash \Phi \rightarrow \Psi} .
\]

\end{proof}

\begin{Lemm} \label{Lemm-QN-Soundness}
Для любой $\Phi \in \mathrm{Form}_{\sigma}$,
\[
\vdash \Phi
\quad \Longrightarrow \quad
\vDash \Phi .
\]
\end{Lemm}

\begin{proof}

Предположим, что $\vdash \Phi$. Тогда $\Phi \in \mathsf{QBK}_{\sigma}$, то есть\ существует конечная последовательность
\[
\Phi_0 , \quad \Phi_1 , \quad \dots , \quad \Phi_n\ =\ \Phi
\]
 $\sigma$-формул такая, что для любого $i \in \left\{ 0, \dots, n \right\}$ верно одно из следующих условий:
\begin{enumerate}

\item[a.] $\Phi_i$ аксиома;

\item[b.] $\Phi_i$ получена из предыдущих $\Phi_j$ и $\Phi_k$ по правилу $\mathtt{MP}$;

\item[c.] $\Phi_i$ получена из предыдущей $\Phi_j$ по правилу $\mathtt{MB}$ или $\mathtt{MD}$;

\item[d.] $\Phi_i$ получена из предыдущей $\Phi_j$ по правилу $\mathtt{BR1}$ или $\mathtt{BR2}$.

\end{enumerate}
Пусть $\mathcal{M}$ --- $\mathsf{QBK}_{\sigma}$-модель. Индукцией по $i$ покажем, что $\mathcal{M}, w \Vdash^+ \lambda(\Phi_i)$ для любых $w \in W$ и  основных $A_w$-подстановок $\lambda$. 
\medskip

Общезначимость аксиом пропозициональной классической логики и аксиом $\mathtt{Q1}$, $\mathtt{Q2}$ проверяется так же, как и в логике предикатов. \vspace{2 mm} \\
Общезначимость пропозициональных аксиом сильного отрицания, за исключением $\mathtt{SN2}$, очевидным образом следует из определения. Проверим $\mathtt{SN2}$: \vspace{2 mm} \\ 
Пусть $\Phi_i={\sim(\Psi\rightarrow\Theta)\leftrightarrow(\Psi\wedge\sim\Theta)}$. Рассмотрим импликацию слева направо $$\mathcal{M},w\Vdash^+\lambda(\sim(\Psi\rightarrow\Theta)\rightarrow(\Psi\wedge\sim\Theta)).$$
Так как $$\lambda(\sim(\Psi\rightarrow\Theta)\rightarrow(\Psi\wedge\sim\Theta))=(\sim(\lambda(\Psi)\rightarrow\lambda(\Theta))\rightarrow(\lambda(\Psi)\wedge\sim\lambda(\Theta)),$$
то, расписывая по определению верифицируемость импликации, имеем $$\mathcal{M},w\Vdash^+\sim(\lambda(\Psi)\rightarrow\lambda(\Theta))\implies\mathcal{M},w\Vdash^+\lambda(\Psi)\wedge\sim\lambda(\Theta).$$
Верифицируя сильное отрицание слева, получаем, что нужно фальсифицировать импликацию. Значит, $$\mathcal{M},w\Vdash^+\lambda(\Psi)~~~\text{и}~~~\mathcal{M},w\Vdash^-\lambda(\Theta).$$
Верифицируя конъюнкцию справа, получаем такое же условие. Общезначимость обратной импликации  следует из проведённого рассуждения.\vspace{2 mm} \\
Также легко проверяется общезначимость остальных аксиом.

\medskip
Если $\Phi_i$ получена из предыдущих формул по правилу $\mathtt{MP}$ или правилам Бернайса, то её общезначимость проверяется снова так же, как и в логике предикатов.

\medskip
Пусть $\Phi_i$ получена из предыдущей $\Phi_j$ по правилу $\mathtt{MB}$, значит, для некоторых $\Psi$ и $\Theta$ имеем $$\Phi_j={\Psi\rightarrow\Theta}~~~\text{и}~~~\Phi_i={\square\Psi\rightarrow\square\Theta}.$$ 
По индукционному предположению $\mathcal{M},w\Vdash^+\lambda(\Psi\rightarrow\Theta)$ для любых $w\in W$ и  основной $A_w$-подстановки $\lambda$. Это означает, что $$\mathcal{M},w\Vdash^+\lambda(\Psi)\implies\mathcal{M},w\Vdash^+\lambda(\Theta)~~~~(*)$$
Нам достаточно показать, что $$\forall w_1(wRw_1\implies\mathcal{M},w_1\Vdash^+\lambda(\Psi))\implies \forall w_1(wRw_1\implies\mathcal{M},w_1\Vdash^+\lambda(\Theta)).$$
Но это следует из того, что $(*)$ верна во всех мирах и при всех основных подстановках. \vspace{2 mm} \\
Случай правила $\mathtt{MD}$ проверяется аналогично.
\end{proof}

\begin{Th}[о корректности] \label{Theorem_Soundness}
Для любых множеств $\Gamma \subseteq \mathrm{Sent}_{\sigma}$ и $\Delta \subseteq \mathrm{Form}_{\sigma}$,
\[
{\Gamma \vdash \Delta}
\quad \Longrightarrow \quad
{\Gamma \vDash \Delta} .
\]
\end{Th}

\begin{proof}
Предположим $\Gamma \vdash \Delta$, то есть\ $\Gamma \vdash \Phi$ для некоторой $\Phi \in \mathrm{Disj} \left( \Delta \right)$. Тогда существует конечное подмножество
 $\Lambda$ множества $\Gamma$ такое, что $\Lambda \vdash \Phi$. Имеем два случая.

\medskip
Пусть $\Lambda = \varnothing$. Тогда $\vDash \Phi$ по лемме 10. Следовательно, $\Gamma \vDash \Delta$.

\medskip
Пусть $\Lambda = \left\{ \Psi_0, \dots, \Psi_n \right\}$. Тогда $\Psi_0 \wedge \dots \wedge \Psi_n \vdash \Phi$, что эквивалентно $\vdash {\Psi_0 \wedge \dots \wedge \Psi_n} \rightarrow \Phi$ по теореме о дедукции. Следовательно, $\vDash {\Psi_0
\wedge \dots \wedge \Psi_n} \rightarrow \Phi$ по лемме 10. Тогда $\Psi_0 \wedge \dots \wedge \Psi_n
\vDash \Phi$ по теореме о семантической дедукции. Следовательно, $\Gamma \vDash \Delta$.
\end{proof}
\section{Теорема о сильной полноте}






Назовём $\sigma$-теорию $\Gamma$ \emph{простой}, если 
\begin{itemize}
    \item $\Gamma \ne \mathrm{Sent}_{\sigma}$;
    \item $\left\{ \Phi \in \mathrm{Sent}_{\sigma} \mid \Gamma \vdash \Phi \right\} \subseteq \Gamma$;
    \item если $\Phi \vee \Psi \in \Gamma$, то $\Phi \in \Gamma$ или $\Psi \in \Gamma$.
\end{itemize}
С простыми теориями связана широко известная лемма.
\begin{Lemm}[Линденбаума]
Пусть $\Gamma$ --- некоторая $\sigma$-теория. Если  $\Gamma\nvdash\bot$,  то $\Gamma$ расширяется до простой $\sigma$-теории.
\end{Lemm}
Назовём простую $\sigma$-теорию $\Gamma$ \emph{насыщенной}, если для любого предложения  ${\exists x}\, \Phi$, содержащегося в $\Gamma$, существует константа $c\in \mathrm{Const}_{\sigma}$ такая, что  $\Phi \left( x/c \right) \in \Gamma$. Такое свойство $\sigma$-теорий называют \emph{экзистенциальным}. \vspace{2 mm} \\
Заметим, что если $\Gamma$ является насыщенной теорией и содержит предложения $\Phi(x/c)$ для любых констант $c$, то $\Gamma$ содержит формулу $\forall x\Phi$.\footnote{Это следует из того, что $\Gamma$ дедуктивно замкнута, содержит предложение $\Phi\vee\neg\Phi$ и обладает экзистенциальным свойством.} Для насыщенных теорий верна стандартная лемма о расширении; смотри, например, \cite{Speranski}.
\begin{Lemm}[о расширении] \label{Extension_Lemma_Strong}
Пусть $\Gamma \subseteq \mathrm{Sent}_{\sigma}$ и $\Delta \subseteq \mathrm{Form}_{\sigma}$ такие, что $\Gamma \nvdash \Delta$. Тогда для любого множества
$S$ мощности $\left| \mathrm{Sent}_{\sigma} \right|$ существует насыщенная $\sigma_{S}$-теория $\Gamma' \supseteq \Gamma$ такая, что $\Gamma' \nvdash \Delta$.
\end{Lemm}

Далее под $(\Gamma)_{\square}$ будем подразумевать множество $\{\Phi~|~\square\Phi\in\Gamma\}$. Для определения канонической модели для $\mathsf{QBK}_{\sigma}$  зафиксируем некоторое множество  $S^{\star}$ мощности $\left| \mathrm{Sent}_{\sigma} \right|$. Назовём $S \subseteq S^{\star}$ \emph{допустимым}, если $\left| S^{\star} \setminus S \right| = \left| S^{\star} \right|$. 
\medskip
Для любого множества $S$ определим
\[
\mathrm{Saturated}_S\ :=\
\text{совокупность всех насыщенных} ~ \sigma_{S} \text{-теорий} .
\]
Свяжем с каждой  $\Gamma \in \mathrm{Saturated}_S$  $\sigma_S$-структуры $(\mathfrak{A}_{\Gamma}^S)^+$ и $(\mathfrak{A}_{\Gamma}^S)^-$ с носителем $\mathrm{Const} \left( \Gamma
\right)$ таким, что всякий константный символ из $\sigma_S$ интерпретируется как он сам, и для любого атомарного  $\sigma_S$-предложения $\Phi$
\[
({\mathfrak{A}_{\Gamma}^S)^+ \Vdash \Phi} \quad :\Longleftrightarrow \quad
{\Phi\ \in\ \Gamma}, ~~~~~~~~~~~~~ ({\mathfrak{A}_{\Gamma}^S)^- \Vdash \Phi} \quad :\Longleftrightarrow \quad
{\sim\Phi\ \in\ \Gamma}.
\]
Заметим, что $\mathrm{Const}(\Gamma)=\mathrm{Const}_{\sigma_{S}}$.\footnote{С одной стороны, $\Gamma$ как $\sigma_{S}$-теория не может содержать формулы с константами не из $\sigma_{S}$. А с другой стороны, так как $\Gamma$ дедуктивно замкнута, в ней содержатся любые выводимые из аксиом формулы. Однако каждая константа из $Const_{\sigma_{S}}$ содержится в какой-нибудь такой формуле. Следовательно, $Const(\Gamma)=Const_{\sigma_{S}}$. } \vspace{2 mm} \\
Обозначим за $\mathfrak{A}_{\Gamma}^+$ и $ \mathfrak{A}_{\Gamma}^-$  $\sigma$-обеднения  $(\mathfrak{A}_{\Gamma}^S)^+ $ и $ (\mathfrak{A}_{\Gamma}^S)^-$ соответственно. Очевидно, каждое \\ $A_{\Gamma}$-предложение имеет вид
\[
{\Phi \left( x_1/\underline{c_1}, \dots, x_n/\underline{c_n} \right)},
\]
где $\left\{ c_1, \dots, c_n \right\} \subseteq \mathrm{Const} \left( \Gamma \right)$. Далее каждое такое предложение  будем отождествлять  с \\ $\sigma_S$-предложением $\Phi
\left( x_1/c_1, \dots, x_n/c_n \right)$.

\medskip
Множеством миров нашей модели положим \[
W^{\mathsf{QBK}}\ :=\
\bigcup {\left\{ \mathrm{Saturated}_S \mid S ~ \text{допустимое подмножество} ~ S^{\star} \right\}} .
\]
Под \emph{канонической шкалой для $\mathsf{QBK}$} будем подразумевать $\mathcal{W}^{\mathsf{QBK}} = \langle W^{\mathsf{QBK}}, R^{\mathsf{QBK}}
\rangle$, где 
\[R^{\mathsf{QBK}}\ :=\
{\left\{ \left( \Gamma_1, \Gamma_2 \right) \in W^{\mathsf{QBK}} \times W^{\mathsf{QBK}} \mid (\Gamma_1)_{\square} \subseteq \Gamma_2  \right\}} .\]
Тогда \emph{каноническая модель для $\mathsf{QBK}$} --- это $\mathcal{M}^{\mathsf{QBK}} = \langle \mathcal{W}^{\mathsf{QBK}}, (\mathscr{A}^{\mathsf{QBK}})^+,(\mathscr{A}^{\mathsf{QBK}})^-
\rangle$, где
\[
{(\mathscr{A}^{\mathsf{QBK}})^+ \left( \Gamma \right)}\ :=\
\mathfrak{A}_{\Gamma}^+,~~~~~~~~{(\mathscr{A}^{\mathsf{QBK}})^- \left( \Gamma \right)}\ :=\
\mathfrak{A}_{\Gamma}^-.
\]
Легко убедиться в том, что это действительно $\mathsf{QBK}$-модель.\footnote{Любая константа из $Const(\Gamma_1)$ содержится в некоторой выводимой из аксиом формуле $\Phi$. По правилу нормализации, $\square\Phi$ также будет выводима из аксиом, а следовательно, лежать в $\Gamma_1$. Тогда по определению $(\Gamma_1)_{\square}$ $\Phi$ будет лежать в $\Gamma_2$. Таким образом, носитель всегда не уменьшается. Интерпретация  констант из $\sigma$ сохраняется тривиально.}
\begin{Lemm}[о канонической модели] 
Для любой $\Gamma \in W^{\mathsf{QBK}}$ и $A_{\Gamma}$-предложения $\Phi$,
\[
\mathcal{M}^{\mathsf{QBK}}, \Gamma \Vdash^+ \Phi
\quad \Longleftrightarrow \quad
\Phi\ \in\ \Gamma,~~~~~~~\mathcal{M}^{\mathsf{QBK}}, \Gamma \Vdash^- \Phi
\quad \Longleftrightarrow \quad
\sim\Phi\ \in\ \Gamma. 
\]
\end{Lemm}

\begin{proof}
Индукция по сложности предложения $\Phi$.

\medskip
Случай атомарного $\Phi$ очевиден.

\medskip
Предположим, что $\Phi=\exists x\Psi$ и докажем сначала случай верифицируемости:\vspace{2 mm} \\
Пусть 
$\mathcal{M}^{\mathsf{QBK}}, \Gamma \Vdash^+ \exists x\Psi$. Значит, существует  $a\in \mathrm{Const}(\Gamma)$ такая, что $\mathcal{M}^{\mathsf{QBK}}, \Gamma \Vdash^+ \Psi(x/a).$
По предположению индукции $\Psi(x/a)\in\Gamma$. И по аксиоме $\mathtt{Q2}$ получаем, что  $\exists x\Psi\in\Gamma$.\vspace{2 mm} \\
Теперь пусть $\exists x \Psi\in\Gamma$. По определению насыщенной теории получаем, что $\Psi(x/c)\in\Gamma$ для некоторой константы из $\mathrm{Const}_{\sigma_S}$. Значит,  $\mathcal{M}^{\mathsf{QBK}}, \Gamma \Vdash^+ \Psi(x/\underline{c})$ ввиду  индукционного предположения.  Следовательно, $\mathcal{M}^{\mathsf{QBK}}, \Gamma \Vdash^+ \exists x\Psi$.
\medskip

Теперь рассмотрим случай фальсифицируемости: \vspace{2 mm} \\
Пусть $\mathcal{M}^{\mathsf{QBK}}, \Gamma \Vdash^- \exists x\Psi$. Тогда для любой $a\in \mathrm{Const}(\Gamma)~~ \mathcal{M}^{\mathsf{QBK}}, \Gamma \Vdash^- \Psi(x/a).$  По индукционному предположению получаем, что $\sim\Psi(x/a)\in\Gamma$ для любой $a\in \mathrm{Const}(\Gamma)$. Но, так как $\mathrm{Const}(\Gamma)=\mathrm{Const}_{\sigma_S}$, вследствие насыщенности имеем $\forall x \sim\Psi\in\Gamma$. Тогда по аксиоме $\mathtt{Q4}$ получаем, что $\sim\exists x\Psi\in\Gamma$.\vspace{2 mm} \\
Теперь предположим, что $\sim\exists x\Psi\in\Gamma$. Тогда по аксиоме $\mathtt{Q4}\,$   $\forall x \sim\Psi\in\Gamma$. Далее по аксиоме $\mathtt{Q1}\,$ $\sim\Psi(x/c)\in\Gamma$ для некоторой константы $c\in \mathrm{Const}_{\sigma_S}$. 
Предположим, что $\mathcal{M}^{\mathsf{QBK}}, \Gamma \nVdash^- \exists x\Psi.$
Значит, для любой константы $c\in \mathrm{Const}(\Gamma)~$ $\mathcal{M}^{\mathsf{QBK}}, \Gamma \nVdash^- \Psi(x/\underline{c})$, то есть  $\mathcal{M}^{\mathsf{QBK}}, \Gamma \nVdash^- \Psi(x/c).$  По индукционному предложению получаем, что для любой константы $c\in \mathrm{Const}(\Gamma)$  $\sim\Psi(x/c)\notin\Gamma$. Противоречие.

\medskip
Теперь предположим, что $\Phi=\forall x\Psi$, и рассмотрим сначала случай верифицируемости:\vspace{2 mm} \\
Пусть $\mathcal{M}^{\mathsf{QBK}}, \Gamma \Vdash^+ \forall x\Psi$. Тогда для любой константы $a\in \mathrm{Const}(\Gamma)~~\mathcal{M}^{\mathsf{QBK}}, \Gamma \Vdash^+ \Psi(x/a).$ По индукционному предположению получаем, что $\Psi(x/a)\in\Gamma$ для любой константы $a$ из $\mathrm{Const}(\Gamma)=\mathrm{Const}_{\sigma_S}$. По свойствам насыщенной теории имеем $\forall x\Psi\in\Gamma$.\vspace{2 mm} \\
С другой стороны, пусть $\forall x\Psi\in\Gamma$. Тогда для любой константы $a\in \mathrm{Const}_{\sigma_S}=\mathrm{Const}(\Gamma)$ по аксиоме $\mathtt{Q1}$  получаем, что $\Psi(x/a)\in\Gamma$. Применяя предположение индукции, имеем $\mathcal{M}^{\mathsf{QBK}}, \Gamma \Vdash^+ \Psi(x/\underline{a}).$ Следовательно, $\mathcal{M}^{\mathsf{QBK}}, \Gamma \Vdash^+ \forall x\Psi$.

\medskip
Рассмотрим случай фальсифицируемости:\vspace{2 mm} \\
Пусть $\mathcal{M}^{\mathsf{QBK}}, \Gamma \Vdash^- \forall x\Psi$. Тогда существует такая константа $a$, что  $\mathcal{M}^{\mathsf{QBK}}, \Gamma \Vdash^- \Psi(x/a)$. По индукционному предположению получаем, что  $\sim\Psi(x/a)\in\Gamma$. Тогда по акcиоме $\mathtt{Q1}$ имеем $\exists x\sim\Psi\in\Gamma$. Наконец, по аксиоме $\mathtt{Q3}$ получаем, что $\sim\forall x\Psi\in\Gamma$.\vspace{2 mm} \\
 Предположим, что $\sim\forall x\Psi\in\Gamma$. Тогда по аксиоме $\mathtt{Q3}$ имеем  $\exists x\sim\Psi\in\Gamma$. Так как $\Gamma$ является насыщенной $\sigma_S$-теорией, получаем, что $\sim\Psi(x/c)\in\Gamma$ для некоторой константы из $\mathrm{Const}_{\sigma_S}$. Пусть $ \mathcal{M}^{\mathsf{QBK}}, \Gamma \nVdash^- \forall x\Psi.$ Значит, для любой константы $a\in \mathrm{Const}(\Gamma)~$ $\mathcal{M}^{\mathsf{QBK}}, \Gamma \nVdash^- \Psi(x/\underline{a})$. Но тогда  $\mathcal{M}^{\mathsf{QBK}}, \Gamma \nVdash^- \Psi(x/a)$
для  тех же констант. По предположению индукции заключаем, что $\sim\Psi(x/a)\notin\Gamma$ для любой константы $a\in \mathrm{Const}(\Gamma)=\mathrm{Const}_{\sigma_S}$. Противоречие.

\medskip
Остальные случаи разбираются так же, как и в логике $\mathsf{BK}$.

\end{proof}

\begin{Th}[о сильной полноте] 
Для любых множеств $\Gamma \subseteq \mathrm{Sent}_{\sigma}$ и $\Delta \subseteq \mathrm{Form}_{\sigma}$,
\[
{\Gamma \vdash \Delta}
\quad \Longleftrightarrow \quad
{\Gamma \vDash \Delta} .
\]
\end{Th}

\begin{proof}
\fbox{$\Longrightarrow$}~\,Это просто теорема о корректности.

\medskip
\fbox{$\Longleftarrow$}~\,Предположим, что $\Gamma \nvdash \Delta$. Зафиксируем некоторое допустимое множество $S \subseteq S^{\star}$ мощности $\aleph_0$ (то есть $\left| S \right|
= \left| \mathrm{Var} \right|$). Пусть $\lambda$ --- биекция между $\mathrm{Var}$ и $\left\{ \underline{s} \mid s \in S \right\}$.
Рассмотрим
\[
{\Delta'}\ :=\
{\left\{ \lambda \Psi \mid \Psi \in \Delta \right\}} .
\]
Но тогда, поскольку новых констант $\underline{s}$ нет ни в $\Gamma$, ни в $\Delta$, получаем, что $\Gamma \nvdash \Delta'$. По лемме о расширении существует  $\Gamma'
\in W^{\mathsf{QBK}}$ такая, что $\Gamma \subseteq \Gamma'$ и $\Gamma' \nvdash \Delta'$. Очевидно, $\lambda$ можно рассматривать как основную
$A_{\Gamma'}$-подстановку. Тогда по лемме 14~ $\mathcal{M}^{\mathsf{QBK}}, \Gamma' \Vdash \Phi$ для любых $\Phi \in \Gamma$,
в то время как $\mathcal{M}^{\mathsf{QBK}}, \Gamma' \nVdash \lambda \Psi$ для всех $\Psi \in \Delta$. Следовательно, $\Gamma \nvDash \Delta$.
\end{proof}
\vspace{7 mm}
\section{Некоторые естественные расширения}
Теперь рассмотрим некоторые естественные расширения $\mathsf{QBK}$, для которых можно установить сильную полноту при помощи $\mathcal{M}^{\mathsf{QBK}}$. Такие расширения можно разделить на две группы: 
\begin{itemize}
    \item[\romannumeral 1] обогащенные новыми аксиомами;
    \item[\romannumeral 2] полученные путем наложения некоторых ограничений на отношение достижимости в шкалах.
\end{itemize} 
\medskip

Для доказательства сильной полноты расширений первой группы в каждом случае необходимо модернизировать понятие модели, проверить корректность новой аксиомы относительно нового понятия модели и удостовериться, что каноническая $\mathsf{QBK}$-модель останется моделью в новом определении. Тогда полнота обогащеной системы будет получатся из доказанных выше леммы о канонической модели и теоремы о сильной полноте для $\mathsf{QBK}$.
\medskip

Итак, обозначим за ${\mathsf{QBK}}^{\circ}_{\sigma}$ наименьшее множество $\sigma$-формул, полученное из $\mathsf{QBK}_{\sigma}$ путём добавления схемы аксиом ${\Phi}\vee{\sim\Phi}$ и замкнутое относительно правил вывода $\mathtt{MP}$, $\mathtt{MB}$, $\mathtt{MD}$, $\mathtt{BR1}$ и $\mathtt{BR2}$. Соответствующее отношение выводимости будем обозначать~ $\vdash_{{\mathsf{QBK}}^{\circ}}$.
\medskip

Назовём упорядоченную тройку $\mathcal{M}=\langle\mathcal{W},\mathscr{A}^+,\mathscr{A}^-\rangle~$ \emph{${\mathsf{QBK}}^{\circ}_{\sigma}$-моделью}, если $\mathcal{M}$ --- это $\mathsf{QBK}_{\sigma}$-модель и для любого мира $w\in W$, и любого атомарного предложения $\Phi$ верно   $\mathfrak{A}^+_w\Vdash\Phi$  или $\mathfrak{A}^-_w\Vdash\Phi$ (или оба сразу). Соответствующее отношение семантического следования  будем обозначать~ $\vDash_{{\mathsf{QBK}}^{\circ}}$.

\begin{Th}[о сильной полноте ${\mathsf{QBK}}^{\circ}$]
Для любых множеств $\Gamma \subseteq \mathrm{Sent}_{\sigma}$ и $\Delta \subseteq \mathrm{Form}_{\sigma}$,
\[
{\Gamma \vdash_{{\mathsf{QBK}}^{\circ}} \Delta}
\quad \Longleftrightarrow \quad
{\Gamma \vDash_{{\mathsf{QBK}}^{\circ}} \Delta} .
\]
\end{Th}
\begin{proof}
Корректность добавленной аксиомы относительно таких моделей проверяется индукцией по построению формулы $\Phi$ с тривиальной базой и простым  шагом индукции. \vspace{2 mm} \\
Проверим, что каноническая $\mathsf{QBK}$-модель будет моделью в смысле ${\mathsf{QBK}}^{\circ}$. Пусть $\Gamma$ --- произвольный мир канонической модели и $\Phi$ --- атомарное предложение. $\Gamma$, как насыщенная теория, обладает дизъюнктивным свойством и содержит формулу ${\Phi}\vee{\sim\Phi}$. Значит, она содержит или $\Phi$, или $\sim\Phi$, или сразу обе формулы. Тогда по лемме о канонической модели выполняется условие на ${\mathsf{QBK}}^{\circ}$-модель.
\end{proof}
Теперь обозначим за $\mathsf{QB3K}_{\sigma}$ наименьшее множество $\sigma$-формул, полученное из $\mathsf{QBK}_{\sigma}$ путём добавления схемы аксиом $\sim{\Phi}\rightarrow(\Phi\rightarrow\Psi)$ и замкнутое относительно правил вывода $\mathtt{MP}$, $\mathtt{MB}$, $\mathtt{MD}$, $\mathtt{BR1}$ и $\mathtt{BR2}$. Соответствующее отношение выводимости  будем обозначать~ $\vdash_{\mathsf{QB3K}}$.
\medskip

Назовём упорядоченную тройку $\mathcal{M}=\langle\mathcal{W},\mathscr{A}^+,\mathscr{A}^-\rangle~$ \emph{$\mathsf{QB3K}_{\sigma}$-моделью}, если $\mathcal{M}$ --- это $\mathsf{QBK}_{\sigma}$-модель, и для любого мира $w\in W$ не существует такого атомарного предложения $\Phi$, что одновременно выполнено $\mathfrak{A}^+_w\Vdash\Phi$ и $\mathfrak{A}^-_w\Vdash\Phi$. Соответствующее отношение семантического следования  будем обозначать~ $\vDash_{\mathsf{QB3K}}$.
\begin{Th}[о сильной полноте $\mathsf{QB3K}$]
Для любых множеств $\Gamma \subseteq \mathrm{Sent}_{\sigma}$ и $\Delta \subseteq \mathrm{Form}_{\sigma}$,
\[
{\Gamma \vdash_{\mathsf{QB3K}} \Delta}
\quad \Longleftrightarrow \quad
{\Gamma \vDash_{\mathsf{QB3K}} \Delta} .
\]
\end{Th}
\begin{proof}
Корректность добавленной аксиомы относительно таких моделей снова проверяется индукцией по построению формулы $\Phi$. \vspace{2 mm} \\
Проверим, что каноническая $\mathsf{QBK}$-модель будет моделью в смысле $\mathsf{QB3K}$. Предположим, что существует такая насыщенная теория $\Gamma$ и атомарное предложение $\Phi$, что одновременно выполнено $\mathfrak{A}^+_{\Gamma}\Vdash\Phi$ и $\mathfrak{A}^-_{\Gamma}\Vdash\Phi$. Тогда по лемме о канонической модели $\Phi$ и $\sim{\Phi}$ лежат в $\Gamma$. Так как $\Gamma$ содержит формулу $\sim{\Phi}\rightarrow(\Phi\rightarrow\Psi)$ и дедуктивно замкнута, получаем, что любое предложение $\Psi$ содержится в $\Gamma$. Значит, $\Gamma=\mathrm{Sent}_{\sigma}$. Но это противоречит определению насыщенной теории\footnote{К $QBK$ также можно добавить две рассмотренные выше схемы одновременно, получив множество формул ${QB3K}^{\circ}$. В определении модели нужно потребовать условия на модели $QB3K$ и ${QBK}^{\circ}$ одновременно. Тогда ${QB3K}^{\circ}$ также окажется сильно полна. }.
\end{proof}
Для доказательства сильной полноты расширений второй группы нужно, во-первых, установить, что отношение достижимости $R$ шкалы $\mathcal{W}$ обладает интересующим нас свойством  тогда и только тогда, когда некоторая модальная формула верифицируема на шкале $\mathcal{W}$. Далее необходимо  показать, что если $\mathsf{QBK}$ содержит такую формулу, то отношение достижимости канонической $\mathsf{QBK}$-модели будет обладать нужным нам свойством. Из этих двух фактов легко получается сильная полнота соответствующих расширений. Доказываются эти два факта аналогично случаю логики $\mathsf{K}$; смотри, например, \cite{Neklassik}.\vspace{2 mm} \\
В частности, при требовании рефлексивности (транзитивности) шкал, мы получим раширение, которое обозначим $\mathsf{QBT}$ ($\mathsf{QBK4}$). Если рассматривать шкалы, обладающие свойством предпорядков, то соответствующее расширение обозначим $\mathsf{QBS4}$. Более полный список логик, для которых можно получить полноту таким способом, читатель снова найдёт в работе \cite{Neklassik}. Обозначим расширение, полученное требованием предпорядка от отношения достижимости в шкале и путём добавления схемы аксиом $\sim\Phi\rightarrow(\Phi\rightarrow\Psi)$, за $\mathsf{QB3S4}$. Теорема о сильной полноте для этого расширения получается как комбинация результатов для $\mathsf{QB3K}$ и $\mathsf{QBS4}$. $\mathsf{QB3S4}$ и $\mathsf{QBS4}$ являются кванторными аналогами соответствующих пропозициональных расширений $\mathsf{BK}$, которые изучались в статье $\cite{Odintsov-Wansing}$ и понадобятся нам в разделе 8.

\vspace{7 mm}
\section{Случай константных носителей}
Рассмотрим формулу $$\mathtt{Ba}:=\lozenge\exists x\Phi\rightarrow\exists x\lozenge\Phi,$$ представляющую собой один из  вариантов так называемой формулы Баркана. Заметим, что  обратная  к ней формула выводима в $\mathsf{QBK}$. Действительно, имеем вывод  \begin{hilbert}

$\Phi\rightarrow\exists x\Phi$& $\mathtt{Q2}$ \\
$\lozenge\Phi\rightarrow\lozenge\exists x\Phi$& $\mathtt{MD}$ к 1\\
$\exists x\lozenge\Phi\rightarrow\lozenge\exists x\Phi$& $\mathtt{BR2}$ к 2

\end{hilbert} 

Обозначим за ${\mathsf{QBK}}^{\sharp}_{\sigma}$ наименьшее множество $\sigma$-формул, полученное из $\mathsf{QBK}_{\sigma}$ путём добавления схемы аксиом $\mathtt{Ba}$ и замкнутое относительно правил вывода $\mathtt{MP}$, $\mathtt{MB}$, $\mathtt{MD}$, $\mathtt{BR1}$ и $\mathtt{BR2}$. Соответствующее отношение выводимости  будем обозначать~ $\vdash_{\sharp}$. 
\vspace{2 mm} 

В ряде случаев вместо формулы $\mathtt{Ba}$ рассматривают другой её вариант $$\mathtt{Ba}^{\square}:=\forall x\square\Phi\rightarrow\square\forall x\Phi.$$
Аналогично формуле $\mathtt{Ba}$ можно показать, что обратный вариант $\mathtt{Ba}^{\square}$ выводим в $\mathsf{QBK}$. При определении $\mathsf{QBK}^{\sharp}$ вместо схемы аксиом $\mathtt{Ba}$ можно добавить к $\mathsf{QBK}$ схему аксиом $\mathtt{Ba}^{\square}$. Обозначим такой вариант $\mathsf{QBK}^{\sharp}_{\mathtt{Ba}^{\square}}$. Однако, как показывает следующее утверждение, доказательство которого практически не отличается от случая классической модальной логики, но приводится для наглядности, два этих варианта образуют одно и то же множество формул.
\begin{Prop}
В $\mathsf{QBK}^{\sharp}$ выводима $\mathtt{Ba}^{\square}$, так же как и в $\mathsf{QBK}^{\sharp}_{\mathtt{Ba}^{\square}}$ выводима $\mathtt{Ba}$.
\end{Prop}
\begin{proof}
Докажем выводимость $\mathtt{Ba}^{\square}$, другой случай доказывается аналогично. \begin{hilbert}
$\lozenge\exists x\neg\Phi\rightarrow\exists x\lozenge\neg\Phi$& частный случай формулы $\mathtt{Ba}$ \\
$\neg\exists x\lozenge\neg\Phi\rightarrow\neg\lozenge\exists x\neg\Phi$& контрапозиция 1 \\
$\forall x\neg\lozenge\neg\Phi\rightarrow\square\neg\exists x\neg\Phi$& классическая тавтология и  $\mathtt{M2}$ к 2 \\
$\forall x\square\neg\neg\Phi\rightarrow\square\neg\exists x\neg\Phi$&  $\mathtt{M2}$ и слабое правило замены к 3 \\
$\square\neg\exists x\neg\Phi\leftrightarrow\square\forall x\neg\neg\Phi$& классические тавтологии и  $\mathtt{MB}$ \\
$\forall x\square\neg\neg\Phi\rightarrow\square\forall x\neg\neg\Phi$& транзитивность импликации  из 4 в 5 \\
$\forall x\square\Phi\rightarrow\forall x\square\neg\neg\Phi$& классическая тавтология, $\mathtt{Q4}$ и правила вывода $\mathtt{MB}$ и $\mathtt{BR1}$ \\ 
$\square\forall x\neg\neg\Phi\rightarrow\square\forall x\Phi$& классическая тавтология,  $\mathtt{Q4}$ и правила вывода $\mathtt{MB}$ и $\mathtt{BR1}$ \\
$\forall x\square\Phi\rightarrow\square\forall x\Phi$& транзитивность импликации из 7 в 6 и из 6 в 8
\end{hilbert}
\end{proof}  Далее будем рассматривать только  $\mathsf{QBK}^{\sharp}$ и ограничимся случаем счётных сигнатур.
\medskip

Упорядоченная тройка $\mathcal{M}=\langle\mathcal{W},\mathscr{A}^+,\mathscr{A}^-\rangle$ называется $\mathsf{QBK}^{\sharp}_{\sigma}$-моделью, если она является $\mathsf{QBK}_{\sigma}$-моделью и для любых двух миров $u,v$ носители соответствующих данным мирам $\sigma$-структур равны. Соответствующее отношение семантического следования  будем обозначать~ $\vDash_{\sharp}$.
\begin{Th}[о корректности ${\mathsf{QBK}}^{\sharp}_{\sigma}$ ] \label{Theorem_Soundness_sharp}
Для любых множеств $\Gamma \subseteq \mathrm{Sent}_{\sigma}$ и $\Delta \subseteq \mathrm{Form}_{\sigma}$,
\[
{\Gamma \vdash_{\sharp} \Delta}
\quad \Longrightarrow \quad
{\Gamma \vDash_{\sharp} \Delta} .
\]
\end{Th}

\begin{proof}
Необходимо проверить только общезначимость формулы $\mathtt{Ba}$. Пусть $\mathcal{M}$ и $w$ --- произвольные ${\mathsf{QBK}}^{\sharp}$-модель и мир, тогда $$\mathcal{M},w\Vdash^+\lozenge\exists x\Phi\rightarrow\exists x\lozenge\Phi\iff(\mathcal{M},w\Vdash^+\lozenge\exists x\Phi\implies\mathcal{M},w\Vdash^+\exists x\lozenge\Phi).$$
При этом левая часть импликации в скобках эквивалентна тому, что существует мир $w'$ такой, что $wRw'$, и элемент $a\in A_{w'}$, что формула $\Phi(x/\underline{a})$ верифицируема в мире $w'$. Тому же эквивалентна правая часть импликации в скобках с одним лишь отличием: $a\in A_{w}$. Так как $A_{w'}=A_{w}$ по определению $\mathsf{QBK}^{\sharp}$-модели, получаем, что импликация верна. 
\end{proof}
\medskip

Перейдём к доказательству сильной полноты $\mathsf{QBK}^{\sharp}$ относительно построенной нами семантики. Для этого модифицируем каноническую модель для $\mathsf{QBK}$.
\vspace{2 mm} \\
Как и ранее, зафиксируем некоторое множество $S^{\star}$ мощности $\left| \mathrm{Sent}_{\sigma} \right|$. Сигнатуру $\sigma_{S^{\star}}$ обозначим как $\sigma^{\star}$. Определим $$\mathrm{Saturated}^{\star}:=\text{совокупность всех насыщенных}\,\sigma^{\star}\text{-теорий}.$$
Аналогично с каждой насыщенной $\sigma^{\star}$-теорией свяжем две $\sigma^{\star}$-структуры. Множеством миров нашей модели положим $$W^{\sharp}=\mathrm{Saturated}^{\star}.$$
При этом каноническая шкала получится из канонической шкалы для $\mathsf{QBK}$ путем релятивизации отношения $R^{\mathsf{QBK}}$ на $W^{\sharp}$. Обозначим это отношение  $R^{\sharp}$, а новую каноническую шкалу  $\mathcal{W}^{\sharp}$.
\vspace{2 mm} \\
Тогда \emph{каноническая модель для $\mathsf{QBK}^{\sharp}$} --- это $\mathcal{M}^{\sharp} = \langle \mathcal{W}^{\sharp}, (\mathscr{A}^{\sharp})^+,(\mathscr{A}^{\sharp})^-
\rangle$, где
\[
{(\mathscr{A}^{\sharp})^+ \left( \Gamma \right)}\ :=\
\mathfrak{A}_{\Gamma}^+,~~~~~~~~{(\mathscr{A}^{\sharp})^- \left( \Gamma \right)}\ :=\
\mathfrak{A}_{\Gamma}^-.
\]
Для таким образом построенной модели и отношения $R^{\sharp}$ верна следующая  лемма, доказательство которой аналогично случаю классической модальной логики; смотри, например, \cite{Gabbay-Shehtman}.
\begin{Lemm}
Для любой $\Gamma\in W^{\sharp}$ и $\Phi\in \mathrm{Sent}_{\sigma^{\star}}$ 
\begin{enumerate}
    \item[\romannumeral 1)] $\lozenge\Phi\in\Gamma$ тогда и только тогда, когда существует $\Gamma'\in W^{\sharp}$ такая, что $\Phi\in\Gamma'$ и  $\Gamma_{\square}\subseteq\Gamma'$;
    \item[\romannumeral 2)] $\square\Phi\in\Gamma$ тогда и только тогда, когда для любых $\Gamma'\in W^{\sharp}$ таких, что $\Gamma_{\square}\subseteq\Gamma'$,  $\Phi\in\Gamma'$.
\end{enumerate}
\end{Lemm} 
\begin{proof}
Сначала докажем  необходимость пункта $\romannumeral 1)$. \vspace{2 mm} \\
Занумеруем множество всех $S^{\star}$-предложений  вида $\exists x \Phi(x)$. Получим множество предложений $\{\exists x\Phi_{k}(x)~|~k>0\}.$ Построим последовательность $(\Gamma_{k}~|~k\in\omega)$, состоящую из конечных $\Gamma_{k}\subseteq Sent_{\sigma^{\star}}$, такую, что для любого $k\in\omega$ выполняется условие $$\lozenge\left(\bigwedge\Gamma_{k}\right)\in\Gamma.$$
Для этого положим $\Gamma_0=\{\Phi\}$, а $\Gamma_{k}=\Gamma_{k-1}\cup\{\exists x\Phi_{k}(x)\rightarrow\Phi_{k}(c)\}$ для некоторой константы $c$ из $S^{\star}$ так, чтобы было выполнено условие выше. Покажем, что такая константа действительно существует.
\medskip

По индукционному предположению имеем $\lozenge\left(\bigwedge\Gamma_{k-1}\right)\in\Gamma.$ Также формула $$\square\exists y(\exists x\Phi_{k}(x)\rightarrow\Phi_{k}(y))$$ лежит в $\mathsf{QBK}^{\sharp}$ для некоторой ещё не использованой переменной $y$ как  классическая тавтология, к которой применили правило обобщения. Но тогда, используя формулу $\lozenge\Phi\wedge\square\Psi\rightarrow\lozenge(\Phi\wedge\Psi)$, которая является тавтологией в пропозициональной модальной логике $\mathsf{K}$, и вынося квантор существования наружу, имеем $$\lozenge\exists y\left(\bigwedge\Gamma_{k-1}\wedge(\exists x\Phi_{k}(x)\rightarrow\Phi_{k}(y))\right)\in\Gamma.$$ По формуле $\mathtt{Ba}$ получаем, что $$\exists y\lozenge\left(\bigwedge\Gamma_{k-1}\wedge(\exists x\Phi_{k}(x)\rightarrow\Phi_{k}(y))\right)\in\Gamma.$$ Так как $\Gamma$ является насыщенной теорией, имеем $$\lozenge\left(\bigwedge\Gamma_{k-1}\wedge(\exists x\Phi_{k}(x)\rightarrow\Phi_{k}(c))\right)\in\Gamma~~\text{для некоторой}~c\in S^{\star}.$$

Положим $\Delta:=\bigcup_{k\in\omega}\Gamma_{k}$. Докажем, что $\Gamma_{\square}\cup\Delta\nvdash_{\sharp}\bot$. Действительно, пусть это не так, тогда $\Gamma_{\square}\cup\Gamma_{k}\vdash_{\sharp}\bot$ для некоторого $k$. По теореме о дедукции $\Gamma_{\square}\vdash_{\sharp}(\bigwedge\Gamma_{k})\rightarrow\bot$, то есть $\Gamma_{\square}\vdash_{\sharp}\neg(\bigwedge\Gamma_{k}).$ Но тогда $\Gamma\vdash_{\sharp}\square\neg(\bigwedge\Gamma_{k}).$ По аксиоме $\mathtt{M2}$ получаем $\Gamma\vdash_{\sharp}\neg\lozenge(\bigwedge\Gamma_{k})$, что невозможно, так как $\Gamma\not = \mathrm{Sent}_{\sigma^{\star}}$.

По лемме Линденбаума получаем, что $\Gamma_{\square}\cup\Delta$ расширяется до простой $\sigma^{\star}$-теории $\Gamma'$, причём $\Phi\in\Gamma'$ и $\Gamma_{\square}\subseteq\Gamma'$. Осталось убедиться, что $\Gamma'$ является насыщенной. Предположим, что $\exists x \Psi(x)\in\Gamma'$, но тогда $\exists x\Psi(x)=\exists x\Phi_{k}(x)$ и $\Psi=\Phi_{k}$ для некоторого $k$. Так как $\exists x\Phi_{k}(x)\rightarrow\Phi_{k}(c)$ лежит в $\Gamma'$, то  по дедуктивной замкнутости $\Gamma'$ получаем, что $\Psi(c)\in\Gamma'$ для некоторой константы $c$. 
\medskip

Докажем теперь достаточность пункта $\romannumeral 1)$. Пусть $\lozenge\Phi\notin\Gamma$. Поскольку $\Gamma$ обладает дизъюнктивным свойством и содержит классическую тавтологию $\neg\Psi\vee\Psi$, имеем $\neg\lozenge\Phi\in\Gamma$. По аксиоме $\mathtt{M2}~$ $\square\neg\Phi\in\Gamma$. Но тогда для любых $\Gamma'\in W^{\sharp}$, если $\Gamma_{\square}\subseteq\Gamma'$, то $\Phi\notin\Gamma'$.
\medskip

Необходимость пункта $\romannumeral 2)$ очевидна, поэтому докажем его достаточность. Пусть $\square\Phi\notin\Gamma$. Тогда, по аналогичным причинам, $\neg\square\Phi\in\Gamma$. По аксиоме $\mathtt{M1}$ получаем, что $\lozenge\neg\Phi\in\Gamma$. По пункту $\romannumeral 1)$ существует $\Gamma'\in W^{\sharp}$ такая, что $\Gamma_{\square}\subseteq\Gamma'$ и $\neg\Phi\in\Gamma'$, то есть  $\Phi\notin\Gamma'$.

\end{proof}
Адаптируем под наши текущие нужды лемму о расширении.
\begin{Lemm}
Пусть $S\subseteq S^{\star}$ допустимо. Предположим, что $\Gamma\subseteq \mathrm{Sent}_{\sigma}$ и $\Delta\subseteq \mathrm{Form}_{\sigma}$ таковы, что $\Gamma\not\vdash\Delta$. Тогда существует насыщенная $\sigma^{\star}$-теория $\Gamma'$ такая, что $\Gamma\subseteq\Gamma'$ и $\Gamma'\not\vdash\Delta$.  
\end{Lemm}
\begin{proof}
Следствие леммы 13 в случае, когда $\sigma=\sigma_{S}$ и $S=S^{\star}$.
\end{proof}
Теперь мы готовы доказать основную семантическую лемму и, следовательно, теорему о сильной полноте.
\begin{Lemm}[о канонической модели] 
Для любой $\Gamma \in W^{\sharp}$ и $A_{\Gamma}$-предложения $\Phi$,
\[
\mathcal{M}^{\sharp}, \Gamma \Vdash^+ \Phi
\quad \Longleftrightarrow \quad
\Phi\ \in\ \Gamma,~~~~~~~\mathcal{M}^{\sharp}, \Gamma \Vdash^- \Phi
\quad \Longleftrightarrow \quad
\sim\Phi\ \in\ \Gamma. 
\]
\end{Lemm}
\begin{proof}
Индукция по сложности предложения $\Phi$. \vspace{2 mm} \\
Нас интересуют только те случаи, где в доказательстве леммы о канонической модели для $\mathsf{QBK}$ применяется стандартная лемма о расширении. Это случаи $\square$ и $\lozenge$. \vspace{2 mm} \\
Верифицируемость этих двух случаев устанавливается напрямую из леммы 20 при помощи индукционного предположения.
\medskip

Рассмотрим случай фальсифицируемости $\square$. По определению и индукционному предположению $$\mathcal{M}^{\sharp}, \Gamma \Vdash^- \square\Phi\iff\exists\Delta\,(\Gamma_{\square}\subseteq\Delta ~\text{и}~\sim\Phi\in\Delta).$$ По лемме 20 это эквивалентно тому, что $\lozenge\sim\Phi\in\Gamma$. Но по аксиоме  $\mathtt{M3}\, $ это эквивалентно $\sim\square\Phi\in\Gamma$.
\medskip

Рассмотрим случай фальсифицируемости $\lozenge$. По определению и индукционному предположению $$\mathcal{M}^{\sharp}, \Gamma \Vdash^- \lozenge\Phi\iff\forall\Delta\,(\Gamma_{\square}\subseteq\Delta \implies\sim\Phi\in\Delta).$$ По лемме 20 это эквивалентно тому, что $\square\sim\Phi\in\Gamma$. Но по аксиоме $\mathtt{M4}\,$ это эквивалентно $\sim\lozenge\Phi\in\Gamma$.
\end{proof}
\begin{Th}[о сильной полноте $\mathsf{QBK}^{\sharp}$] 
Для любых множеств $\Gamma \subseteq \mathrm{Sent}_{\sigma}$ и $\Delta \subseteq \mathrm{Form}_{\sigma}$,
\[
{\Gamma \vdash_{\sharp} \Delta}
\quad \Longleftrightarrow \quad
{\Gamma \vDash_{\sharp} \Delta} .
\]
\end{Th}

\begin{proof}
С точностью до замены леммы 14 на лемму 22 и леммы о расширении на лемму 21 аналогично доказательству теоремы о сильной полноте для $\mathsf{QBK}$.
\end{proof}

Стоит отметить, что, как и в случае $\mathsf{QBK}$, можно также рассматривать естественные расширения $\mathsf{QBK}^{\sharp}$. Вводятся они аналогично расширениям $\mathsf{QBK}$ из секции 6, и для них также остаются верны соответствующие теоремы о полноте.
\vspace{7 mm}

\section{Точное вложение кванторных логик Нельсона}
Для предикатных логик Нельсона $\mathsf{QN4}^{\bot}$ и $\mathsf{QN3}$ известен полезный результат о негативных нормальных формах, аналогичный теореме 8 с точностью до замены сильной эквивалентности на слабую, а именно
\begin{Th}
 Для любой формулы $\Phi$ языка логик Нельсона существует формула $\overline{\Phi}$ в негативной нормальной форме такая, что $$\Phi\leftrightarrow\overline{\Phi}\in\mathsf{QN4}^{\bot}(\mathsf{QN3}).$$
\end{Th}
Для всякой формулы языка логик Нельсона $\overline{\Phi}$, находящейся в негативной нормальной форме, определим преобразование $\tau$, сопоставляющее ей некоторую формулу в модальном языке $\mathsf{QBK}$, рекурсивно следующим образом:
\begin{align*}
\tau P(t_1,\ldots,t_n)~&=~\square P(t_1,\ldots,t_n); \\
\tau \sim P(t_1,\ldots,t_n)~&=~\sim\lozenge P(t_1,\ldots,t_n); \\
\tau(\Phi\vee\Psi)~&=~\tau\Phi\vee\tau\Psi; \\
\tau(\Phi\wedge\Psi)~&=~\tau\Phi\wedge\tau\Psi; \\
\tau(\Phi\rightarrow\Psi)~&=~\square(\tau\Phi\rightarrow\tau\Psi); \\
\tau\bot~&=~\bot; \\
\tau\exists x\Phi~&=~\exists x\,\tau\Phi; \\
\tau\forall x\Phi~&=~\square\forall x\,\tau\Phi.
\end{align*}

Модели $\mathsf{QN4}^{\bot}$ ($\mathsf{QN3}$) и $\mathsf{QBS4}$ ($\mathsf{QB3S4})$ имеют одну и ту же форму $\mathcal{M}=\langle W,\leqslant,\mathscr{A}^+,\mathscr{A}^-\rangle$, где $\leqslant$ является предпорядком на $W\times W$. Однако на модели $\mathsf{QN4}^{\bot}$ ($\mathsf{QN3}$) накладывается дополнительное ограничение: если $u,v$ являются мирами некоторой $\mathsf{QN4}^{\bot}$-модели ($\mathsf{QN3}$-модели), а $\mathfrak{A}^{\pm}_{u}$ и $\mathfrak{A}^{\pm}_{v}$ соответствующие этим мирам $\sigma$-структуры, то $$u\leqslant v\implies P^{\mathfrak{A}^{\pm}_{u}}\subseteq P^{\mathfrak{A}^{\pm}_{v}}~\text{для любого}~P\in \mathrm{Pred}_{\sigma}.$$
Отметим, что в  случае $\mathsf{QN3}$ мы также имеем ещё одно ограничение: $P^{\mathfrak{A}^{+}_{u}}\cap P^{\mathfrak{A}^{-}_{u}}=\emptyset$ для любых миров $u$.
\medskip

Здесь и далее условимся обозначать за $\overline{\Phi}$ формулу языка логик Нельсона, находящуюся в негативной нормальной форме, а за $\tau(\overline{\Phi})$ результат применения к ней преобразования $\tau$. Теоремы о сильной полноте относительно подходящих семантик Крипке для $\mathsf{QN4}^{\bot}$ и $\mathsf{QN3}$ известны и доказаны в работах \cite{fullN4} и \cite{Gurevich}. \vspace{2 mm} \\
Пусть $\mathcal{M}=\langle W,\leqslant,\mathscr{A}^+,\mathscr{A}^-\rangle$ является $\mathsf{QBS4}$-моделью. Определим $\mathcal{M}'=\langle W,\leqslant,(\mathscr{A}^+)',(\mathscr{A}^-)'\rangle$ следующим образом: 
\begin{itemize}
    \item[\romannumeral 1)] $\mathcal{M}',w\Vdash^+P(\underline{c_1},\ldots,\underline{c_n})\iff\mathcal{M},w\Vdash^+\square P(\underline{c_1},\ldots,\underline{c_n})$;
    \item[\romannumeral 2)] $\mathcal{M}',w\Vdash^-P(\underline{c_1},\ldots,\underline{c_n})\iff\mathcal{M},w\Vdash^-\lozenge P(\underline{c_1},\ldots,\underline{c_n})$.
\end{itemize}
 Очевидно, что построенная таким образом $\mathcal{M}'$ является $\mathsf{QN4}^{\bot}$-моделью. Если бы изначально $\mathcal{M}$ была $\mathsf{QB3S4}$-моделью, то, аналогично, при таком определении $\mathcal{M}'$ получилась бы $\mathsf{QN3}$-моделью.
 \begin{Lemm}
 Пусть $\mathcal{M}=\langle W,\leqslant,\mathscr{A}^+,\mathscr{A}^-\rangle$ является $\mathsf{QBS4}$-моделью. Тогда для любой $\overline{\Phi}$ и мира $w$ $$\mathcal{M}',w\Vdash^+\overline{\Phi}\iff\mathcal{M},w\Vdash^+\tau(\overline{\Phi}).$$ 
 \end{Lemm}
 \begin{proof}
 Индукция по сложности формулы $\overline{\Phi}$.
 \medskip
 
 Случай атомарного $\overline{\Phi}$ очевиден. Случай $\overline{\Phi}={\sim P(\underline{c_1},\ldots,\underline{c_n})}$ легко следует из определения $\mathcal{M}'$: $$\mathcal{M}',w\Vdash^+\sim P(\underline{c_1},\ldots,\underline{c_n})\iff\mathcal{M}',w\Vdash^-P(\underline{c_1},\ldots,\underline{c_n})\iff\mathcal{M},w\Vdash^-\lozenge P(\underline{c_1},\ldots,\underline{c_n})$$ $$\iff\mathcal{M},w\Vdash^+\sim\lozenge P(\underline{c_1},\ldots,\underline{c_n})=\tau(\overline{\Phi}).$$
 Остальные случаи получаются схожим образом. Для примера рассмотрим  $\overline{\Phi}=\forall x\overline{\Psi}.$ Предположим, что утверждение леммы верно для $\overline{\Psi}$. Имеем следующую цепочку эквивалентностей:  $$\mathcal{M}',w\Vdash^+\forall x\overline{\Psi}\iff\forall\,w'\geqslant w\,\forall\,a\in A_{w'}\,\mathcal{M}',w'\Vdash^+\overline{\Psi}(x/\underline{a})$$ $$\iff\forall\,w'\geqslant w\,\forall\,a\in A_{w'}\,\mathcal{M},w'\Vdash^+\tau(\overline{\Psi})(x/\underline{a})\iff\forall\,w'\geqslant w\,\mathcal{M},w'\Vdash^+\forall x\,\tau(\overline{\Psi})$$ $$\iff\mathcal{M},w'\Vdash^+\square\forall x\,\tau(\overline{\Psi})=\tau(\overline{\Phi}).$$
 \end{proof}
 \begin{Th}[о точном вложении предикатных логик Нельсона]
   Формула $\overline{\Phi}$ принадлежит  $\mathsf{QN4}^{\bot}$ ($\mathsf{QN3}$) тогда и только тогда, когда $\tau(\overline{\Phi})$ принадлежит $\mathsf{QBS4}$ ($\mathsf{QB3S4}$).
 \end{Th}
 \begin{proof}
 Предположим, что $\tau(\overline{\Phi})\notin \mathsf{QBS4}$. По теореме о сильной полноте получаем, что $\mathcal{M},w\nVdash^+\tau(\overline{\Phi})$ для некоторой $\mathsf{QBS4}$-модели $\mathcal{M}$ и мира $w$. Но тогда по лемме 24 имеем $\mathcal{M}',w\nVdash^+\overline{\Phi}$. Значит, $\overline{\Phi}\notin\mathsf{QN4}^{\bot}$. 
 \medskip
 
 Теперь пусть $\overline{\Phi}\notin\mathsf{QN4}^{\bot}$. По теореме о сильной полноте имеем $\mathcal{M},w\nVdash^+\overline{\Phi}$ для некоторой $\mathsf{QN4}^{\bot}$-модели и мира $w$. Очевидно, что $\mathcal{M}$ является и $\mathsf{QBS4}$-моделью тоже. Рассматривая $\mathcal{M}$ как $\mathsf{QBS4}$-модель, легко понять, что $\mathcal{M}'=\mathcal{M}$. Тогда по лемме 24 мы тут же получаем, что $\mathcal{M},w\nVdash^+\tau(\overline{\Phi})$, то есть $\tau(\overline{\Phi})\notin \mathsf{QBS4}$.
 \medskip
 
 Случай $\mathsf{QN3}$ и $\mathsf{QB3S4}$ получается аналогично.
 \end{proof}
Из теоремы выше непосредственно следует, что формула $\Phi$ принадлежит $\mathsf{QN4}^{\bot}$ ($\mathsf{QN3}$) тогда и только тогда, когда $\tau(\overline{\Phi})$ принадлежит $\mathsf{QBS4}$ ($\mathsf{QBS4}$). То есть логики Нельсона $\mathsf{QN4}^{\bot}$ ($\mathsf{QN3}$) точно вкладываются в подходящие расширения $\mathsf{QBK}$.

\begin{Cor}[о вложении семантического следования] Для любого множества формул языка логик Нельсона $\Gamma$ и формулы $\Phi$ верно $$\Gamma\vDash_{\star}\Phi\iff\tau(\overline{\Gamma})\vDash_{\widetilde{\star}}\tau(\overline{\Phi}),$$ 
где $\tau(\overline{\Gamma})=\{\tau(\overline{\Psi})~|~\Psi\in\Gamma~\text{и}~\Psi\leftrightarrow\overline{\Psi}\in\mathsf{QN4}^{\bot}(\mathsf{QN3})\}$ и либо $\star=\mathsf{QN4}^{\bot}$ и $\widetilde{\star}=\mathsf{QBS4}$, либо $\star=\mathsf{QN3}$ и $\widetilde{\star}=\mathsf{QB3S4}$.
\end{Cor}
\begin{proof}
Докажем  случай $\mathsf{QN4}^{\bot}$, другой случай доказывается аналогично.
\medskip

\fbox{$\Longrightarrow$}~\, Предположим, что $\Gamma\vDash_{\mathsf{QN4}^{\bot}}\Phi$. Пусть $\mathcal{M}$ и $x$ --- произвольные $\mathsf{QBS4}$-модель и мир. Предположим, что для любой формулы $\tau(\overline{\Psi})\in\tau(\overline{\Gamma})$ верно $\mathcal{M},x\Vdash^+\tau(\overline{\Psi})$. Тогда по лемме 25 $\mathcal{M}',x\Vdash^+\overline{\Psi}$. Так как $\overline{\Psi}$ является негативной нормальной формой для $\Psi$, имеем $\mathcal{M}',x\Vdash^+\Psi$. Значит, для любой формулы $\Psi\in\Gamma$ получаем, что $\mathcal{M}',x\Vdash^+\Psi$. И так как $\Gamma\vDash_{\mathsf{QN4}^{\bot}}\Phi$, имеем $\mathcal{M}',x\Vdash^+\Phi$. Это влечёт $\mathcal{M}',x\Vdash^+\overline{\Phi}$ и тогда снова по лемме 25 $\mathcal{M},x\Vdash^+\tau(\overline{\Phi})$. Так как модель и мир были произвольны, получили заключение импликации.

\fbox{$\Longleftarrow$}~\, Предположим, что $\tau(\overline{\Gamma})\vDash_{\mathsf{QBS4}}\tau(\overline{\Phi})$. Пусть $\mathcal{M}$ и $x$ --- произвольные $\mathsf{QN4}^{\bot}$-модель и мир. Предположим, что для любой формулы $\Psi\in\Gamma$ верно $\mathcal{M},x\Vdash^+\Psi$. Тогда $\mathcal{M},x\Vdash^+\overline{\Psi}$. Очевидно, что $\mathcal{M}$ можно рассматривать как $\mathsf{QBS4}$-модель, и тогда $\mathcal{M}=\mathcal{M}'$. С учётом этого и леммы 25 получаем, что $\mathcal{M},x\Vdash^+\tau(\overline{\Psi})$. Иначе говоря, мы получили, что для любой формулы $\tau(\overline{\Psi})\in\tau(\overline{\Gamma})$ верно $\mathcal{M},x\Vdash^+\tau(\overline{\Psi})$. И так как $\tau(\overline{\Gamma})\vDash_{\mathsf{QBS4}}\tau(\overline{\Phi})$, имеем $\mathcal{M},x\Vdash^+\tau(\overline{\Phi})$. По лемме 25 это эквивалентно тому, что $\mathcal{M}',x\Vdash^+(\overline{\Phi})$. Отсюда получаем, что $\mathcal{M}',x\Vdash^+\Phi$. Вспоминая, что $\mathcal{M}'=\mathcal{M}$, имеем $\mathcal{M},x\Vdash^+\Phi$. Так как модель и мир были произвольны, получили заключение обратной импликации.
\end{proof}
По модулю теорем о сильной полноте для $\mathsf{QN4}^{\bot}$ и $\mathsf{QN3}$ и для соответствующих им расширений $\mathsf{QBK}$ мы также получаем, что $$\Gamma\vdash_{\star}\Phi\iff\tau(\overline{\Gamma})\vdash_{\widetilde{\star}}\tau(\overline{\Phi}).$$ 
\begin{Rem}
Для доказательства точного вложения можно было построить явную трансляцию $\widetilde{\tau}$, избегая непосредственного применения негативных нормальных форм. Она приведена ниже. Однако при внимательном взгляде на данную трансляцию становится ясно, что внутри неё <<зашит>> алгоритм приведения произвольной формулы языка логик Нельсона к негативной нормальной форме. Также легко понять, что результат трансляции $\widetilde{\tau}$ и алгоритма, описанного выше, будут давать один и тот же результат $\tau(\overline{\Phi})$. 

\begin{align*}
\widetilde{\tau} P(t_1,\ldots,t_n)~&=~\square P(t_1,\ldots,t_n) &\widetilde{\tau} \sim P(t_1,\ldots,t_n)~&=~\sim\lozenge P(t_1,\ldots,t_n) \\
\widetilde{\tau}(\Phi\vee\Psi)~&=~\widetilde{\tau}\Phi\vee\widetilde{\tau}\Psi &\widetilde{\tau}\sim(\Phi\vee\Psi)~&=~\widetilde{\tau}\sim\Phi\wedge\widetilde{\tau}\sim\Psi \\
\widetilde{\tau}(\Phi\wedge\Psi)~&=~\widetilde{\tau}\Phi\wedge\widetilde{\tau}\Psi &\widetilde{\tau}\sim(\Phi\wedge\Psi)~&=~\widetilde{\tau}\sim\Phi\vee\widetilde{\tau}\sim\Psi \\
\widetilde{\tau}(\Phi\rightarrow\Psi)~&=~\square(\widetilde{\tau}\Phi\rightarrow\widetilde{\tau}\Psi) &\widetilde{\tau}\sim(\Phi\rightarrow\Psi)~&=~\widetilde{\tau}\Phi\wedge\widetilde{\tau}\sim\Psi \\
\widetilde{\tau}\bot~&=~\bot &\widetilde{\tau}\sim\sim\Phi~&=~\widetilde{\tau}\Phi \\
\widetilde{\tau}\exists x\Phi~&=~\exists x\,\widetilde{\tau}\Phi &\widetilde{\tau}\sim\exists x\Phi~&=~\square\forall x\,\widetilde{\tau}\sim\Phi \\
\widetilde{\tau}\forall x\Phi~&=~\square\forall x\,\widetilde{\tau}\Phi &\widetilde{\tau}\sim\forall x\Phi~&=~\exists x\,\widetilde{\tau}\sim\Phi 
\end{align*}
\end{Rem} 

\begin{Rem}
В работе \cite{Thomason} в разделе Applications, с точностью до замены сильного отрицания в правых частях равенств на классическое и $$\widetilde{\tau} \sim P(t_1,\ldots,t_n)=\square\neg P(t_1,\ldots,t_n),$$  приведена описанная выше трансляция. Обозначим её $\tau'$. Там же утверждается, что она точно вкладывает предикатную логику Нельсона $\mathsf{QN3}+\mathtt{CD}$ в $\mathsf{QS4}+\mathtt{Ba}$. Здесь под $\mathtt{CD}$ имеется в виду следующая схема аксиом: $${\forall x}\, {\left( \Phi \vee \Psi \right)} \rightarrow {\Phi \vee {\forall x}\, \Psi},~~\text{где}~x~\text{не свободен в}~\Phi.$$  Однако это ошибочное утверждение. Докажем, что даже её ограничение на пропозициональную часть не является точным вложением $\mathsf{N3}$ в $\mathsf{S4}$.

Рассмотрим следующую формулу языка пропозициональных логик  Нельсона $$\phi:=(p\rightarrow\sim p)\rightarrow \sim p.$$
Докажем, что она не лежит в $\mathsf{N3}$. Рассмотрим шкалу $W=\langle \{x\},\{(x,x)\}\rangle$ и связаную с ней модель $\mathcal{M}=\langle W,v^+,v^-\rangle$, где  $v^+(p)$ и $v^-(p)$ --- пустые множества, а на остальных переменных оценки определены произвольным образом. Опровергнем формулу $\phi$ в этой модели в мире $x$. Посылка внешней импликации истинна, поскольку в $x$ переменная $p$ не верифицируема и не фальсифицируема. Однако заключение ложно, поскольку в $x$ переменная $p$ не фальсифицируема. Таким образом, вся формула не верифицируема в $x$, следовательно, по теореме о корректности не лежит в  $\mathsf{N3}$.
\medskip

Теперь докажем, что $\tau'(\phi)$ общезначима в  $\mathsf{S4}$. Непосредственно проверяется, что $$\tau'(\phi)=\square(\square(\square p\rightarrow\square\neg p)\rightarrow\square\neg p).$$
Для начала проверим, что $\square(\square p\rightarrow\square\neg p)\rightarrow\square\neg p$ общезначима в $\mathsf{S4}$. Пусть $\mathcal{M}$ и $x$ --- произвольные $\mathsf{S4}$-модель и мир. Необходимо показать, что $$\mathcal{M},x\Vdash\square(\square p\rightarrow\square\neg p)\implies\mathcal{M},x\Vdash\square\neg p.$$
Предположим, что заключение ложно, то есть $\mathcal{M},x\nVdash\square\neg p.$ Значит, существует мир $z$, достижимый из $x$, такой, что $z\in v(p)$. Докажем, что в таком случае посылка тоже ложна $$\mathcal{M},x\nVdash\square(\square p\rightarrow\square\neg p)\iff\exists x'\geqslant x~(\mathcal{M},x'\Vdash\square p~\wedge~\mathcal{M},x'\nVdash\square\neg p).$$
Возьмём $x'=z$. Тогда остаётся показать, что $\mathcal{M},z\Vdash\square p.$ Но это очевидно, поскольку в $\mathsf{S4}$ $v(p)$  является конусом, $z\in v(p)$ и отношение достижимости есть предпорядок.  \vspace{2 mm} \\ 
Пользуясь теоремой о сильной полноте для $\mathsf{S4}$, получаем, что $\square(\square p\rightarrow\square\neg p)\rightarrow\square\neg p\in \mathsf{S4}$. По правилу нормализации $\tau'(\phi)$ лежит в $\mathsf{S4}$, но тогда $\tau'(\phi)$ является общезначимой.
\medskip

Таким образом, мы нашли формулу языка пропозициональных логик Нельсона, которая опровержима в $\mathsf{N3}$, однако её трансляция общезначима в $\mathsf{S4}$. Этот факт опровергает утверждение о точном вложении $\mathsf{N3}$ в $\mathsf{S4}$ посредством трансляции $\tau'$.
\end{Rem}
\vspace{7 mm}

\end{document}